\documentclass[11pt]{amsart}
\usepackage{amsmath}
\usepackage{amsxtra}
\usepackage{amscd}
\usepackage{amsthm}
\usepackage{amsfonts}
\usepackage{amssymb}
\usepackage{eucal}

\newtheorem{cor}{Corollary}
\newtheorem{lem}{Lemma}
\newtheorem{prop}{Proposition}

\theoremstyle{definition}
\newtheorem{definition}{Definition}
\newtheorem{thm}{Theorem}
\newtheorem{conj}{Conjecture}

\theoremstyle{remark}
\newtheorem{rem}{Remark}[section]

\numberwithin{equation}{section}

\begin{document}

\newcommand{\thmref}[1]{Theorem~\ref{#1}}
\newcommand{\secref}[1]{Sect.~\ref{#1}}
\newcommand{\lemref}[1]{Lemma~\ref{#1}}
\newcommand{\propref}[1]{Proposition~\ref{#1}}
\newcommand{\corref}[1]{Corollary~\ref{#1}}
\newcommand{\remref}[1]{Remark~\ref{#1}}
\newcommand{\conjref}[1]{Conjecture~\ref{#1}}
\newcommand{\nc}{\newcommand}
\nc{\on}{\operatorname}
\nc{\ch}{\mbox{ch}}
\nc{\Z}{{\mathbb Z}}
\nc{\C}{{\mathbb C}}
\nc{\pone}{{\mathbb C}{\mathbb P}^1}
\nc{\pa}{\partial}
\nc{\F}{{\mathcal F}}
\nc{\arr}{\rightarrow}
\nc{\larr}{\longrightarrow}
\nc{\al}{\alpha}
\nc{\ri}{\rangle}
\nc{\lef}{\langle}
\nc{\W}{{\mathcal W}}
\nc{\la}{\lambda}
\nc{\ep}{\epsilon}
\nc{\su}{\widehat{{\mathfrak sl}}_2}
\nc{\sw}{{\mathfrak s}{\mathfrak l}}
\nc{\g}{{\mathfrak g}}
\nc{\h}{{\mathfrak h}}
\nc{\n}{{\mathfrak n}}
\nc{\N}{\widehat{\n}}
\nc{\G}{\widehat{\g}}
\nc{\De}{\Delta_+}
\nc{\gt}{\widetilde{\g}}
\nc{\Ga}{\Gamma}
\nc{\one}{{\mathbf 1}}
\nc{\z}{{\mathfrak Z}}
\nc{\zz}{{\mathcal Z}}
\nc{\Hh}{{\mathcal H}_\beta}
\nc{\qp}{q^{\frac{k}{2}}}
\nc{\qm}{q^{-\frac{k}{2}}}
\nc{\La}{\Lambda}
\nc{\wt}{\widetilde}
\nc{\wh}{\widehat}
\nc{\qn}{\frac{[m]_q^2}{[2m]_q}}
\nc{\cri}{_{\on{cr}}}
\nc{\kk}{h^\vee}
\nc{\sun}{\widehat{\sw}_N}
\nc{\hh}{{\mathbf H}_{q,t}}
\nc{\HH}{{\mathcal H}_{q,t}}
\nc{\hhh}{{\mathcal H}_{q,1}}
\nc{\ca}{\wt{{\mathcal A}}_{h,k}(\sw_2)}
\nc{\si}{\sigma}
\nc{\gl}{\widehat{{\mathfrak g}{\mathfrak l}}_2}
\nc{\el}{\ell}
\nc{\s}{T}
\nc{\bi}{\bibitem}
\nc{\om}{\omega}
\nc{\WW}{\W_\beta}
\nc{\scr}{{\mathbf S}}
\nc{\ab}{{\mathbf a}}
\nc{\rr}{r}
\nc{\ol}{\overline}
\nc{\con}{qt^{-1} + q^{-1}t}
\nc{\den}{q^{\el-1} t^{-\el+1}+ q^{-\el+1} t^{\el-1}}
\nc{\ds}{\displaystyle}
\nc{\B}{B}
\nc{\A}{A^{(2)}_{2\el}}
\nc{\GG}{{\mathcal G}}
\nc{\UU}{{\mathcal U}}
\nc{\MM}{{\mathcal M}}
\nc{\CC}{{\mathcal C}}
\nc{\GL}{^L\G}
\nc{\gL}{^L\g}
\nc{\dzz}{\frac{dz}{z}}
\nc{\Res}{\on{Res}}
\nc{\rep}{{\mathcal R}ep \;}
\nc{\uqg}{U_q \G}
\nc{\uqgg}{U_q \g}
\nc{\mc}{\mathcal}
\nc{\Cal}{\mathcal}
\nc{\calp}{{\Cal P}}
\nc{\bp}{{\mathbf P}}
\nc{\bq}{{\mathbf Q}}
\nc{\bb}{{\mathfrak b}}
\nc{\uqb}{U_q \bb_-}
\nc{\uqn}{U_q \wt{{\mathfrak n}}}
\nc{\uqh}{U_q \wt{{\mathfrak h}}}
\nc{\uqhh}{U_q \wh{{\mathfrak h}}}
\nc{\uqnn}{U_q \wh{{\mathfrak n}}}
\nc{\ot}{\otimes}
\nc{\R}{{\mc R}}
\nc{\uqbb}{U_q \wt{\g}}
\nc{\yy}{{\mc Y}}
\nc{\uqsl}{U_q \widehat{\sw}_2}
\nc{\ga}{\gamma}
\nc{\Ab}{{\mathbf A}}
\nc{\Yb}{{\mathbf Y}}
\nc{\yb}{{\mathbf y}}
\nc{\uh}{U \wt{{\mathfrak h}}}
\nc{\uhh}{U \wh{{\mathfrak h}}}
\newcommand{\scs}{\scriptstyle}
\nc{\us}{\underset}
\nc{\opl}{\oplus}
\nc{\yyy}{\wh{\yy}}
\nc{\ovl}{\overline}
\nc{\beq}{\begin{equation}}
\nc{\Fq}{{\mathcal F}}
\nc{\Mq}{{\mathcal M}}
\nc{\Rep}{\on{Rep}}

\title[$q$--characters and ${\mathcal W}$--algebras]{The
$q$--characters of representations of quantum affine algebras and
deformations of ${\mathcal W}$--algebras}

\author[Edward Frenkel]{Edward Frenkel$^*$}\thanks{$^*$Packard Fellow}

\author{Nicolai Reshetikhin}

\address{Department of Mathematics, University of California, Berkeley, CA
94720, USA}

\date{October 1998}

\begin{abstract}
We propose the notion of $q$--characters for finite-dimensional
representations of quantum affine algebras. It is motivated by our
theory of deformed $\W$--algebras.
\end{abstract}

\maketitle

\section{Introduction}

Let $\g$ be a simple Lie algebra, $\G$ be the corresponding
non-twisted affine Kac-Moody algebra, and $\uqg$ be its quantized
universal enveloping algebra (in this paper, $q$ is assumed to be
generic). Consider the category $\rep \uqg$, whose objects are the
finite-dimensional representations of $\uqg$, and morphisms are
homomorphisms of $\uqg$--modules. Since $\uqg$ is a Hopf algebra,
$\rep \uqg$ is a monoidal tensor category.

An interesting problem is to describe the irreducible objects of
${\mathcal R}ep \; \uqg$. A complete answer is known in the case when
$\g=\sw_2$ \cite{CP2} (it is recalled in \secref{casesl2}). For $\g$
other than $\sw_2$ the picture is less clear (see, e.g.,
\cite{CP,CP3,CP4,AK,GV,V}). In contrast, when $q=1$, the analogous
problem of describing irreducible finite-dimensional representations
of $\G$ has a simple and elegant solution. Consider the ``evaluation
homomorphism'' $\phi_a: \G \arr \g$ corresponding to evaluating a
function on $\C^\times$ at a point $a \in \C^\times$. For an
irreducible $\g$--module $V_\la$ with highest weight $\la$, let
$V_\la(a)$ be its pull-back under $\phi_a$ to an irreducible
representation of $\G$. Then $V_{\la_1}(a_1) \otimes \ldots \otimes
V_{\la_n}(a_n)$ is irreducible if $a_i \neq a_j, \forall i\neq j$, and
these are all irreducible finite-dimensional representations of $\G$
up to an isomorphism \cite{CP1}. It is also easy to decompose tensor
products of such representations.

The first major difference, which makes the description of irreducible
representations difficult in the quantum case is that the evaluation
homomorphisms $\phi_z: \G \arr \g$ can not be lifted to homomorphisms
$\uqg \arr \uqgg$ if $\g$ is not $\sw_N$. Nevertheless, V.~Chari and
A.~Pressley have shown \cite{CP,CP3} that for each
$i=1,\ldots,\el=\on{rk} \g$, and $z \in \C^\times$, there exists a
unique irreducible representation $V_{\omega_i}(a)$ of $\uqg$ with
highest weight $\omega_i$, when restricted to $\uqgg \subset
\uqg$. This representation plays the role of $V_{\omega_i}(a)$ in the
case $q=1$, though in general it is bigger than $V_{\omega_i}$ when
restricted to $U_q \g$. Furthermore, Chari and Pressley have shown
\cite{CP3} that any irreducible representation occurs as a subquotient
of the tensor product $V_{\omega_{i_1}}(a_1) \otimes \ldots \otimes
V_{\omega_{i_n}}(a_n)$, where the parameters
$(\omega_{i_1},a_1),\ldots,(\omega_{i_n},a_n)$ are uniquely determined
by this representation up to permutation.

This does provide us with a good parametrization of irreducible
representations, but does not quite answer the question of describing
these representations and their tensor products explicitly. For
instance, the only thing that is known about the tensor product
$V_{\omega_{i}}(a_1) \otimes V_{\omega_{j}}(b)$ in general is that it
is irreducible provided that $a/b$ does not belong to a countable set
(see \cite{KS}), which is unknown in general. This constitutes the
second major difference with the case $q=1$, when this set consists of
a single element, $1$.

In order to gain some insights into the problem, we develop in this
paper a theory of ``characters'' for finite-dimensional
representations of $\uqg$, which we call the $q$--{\em characters}.

Let us recall the situation in the case of finite-dimensional
representations of the Lie algebra $\g$. Such representations can be
integrated to representations of the simply-connected Lie group
$G$. Let $\on{Rep} G$ be the Grothendieck ring of finite-dimensional
representations of $G$. Denote by $T$ the Cartan subgroup of $G$. We
attach to each finite-dimensional $G$--module $V$ its character, the
function $\chi_V: T \arr \C$, defined by $\chi_V(t) =
\on{Tr}_V(t),\forall t \in T$. This way we obtain an injective
homomorphism of commutative algebras $$\chi: \on{Rep} G \arr \Z[T]
\simeq \Z[y_1^{\pm 1},\ldots,y_\el^{\pm 1}],$$ where $y_i$ are the
fundamental weights (the generators of the lattice of homomorphisms $T
\arr \C^\times$).

Let $\on{Rep} \uqg$ be the Grothendieck ring of the category $\rep
\uqg$. Using the universal $R$--matrix of $\uqg$, we will construct an
injective homomorphism
$$\chi_q: \on{Rep} U_q \G \arr \yy = \Z[Y_{i,a_i}^{\pm
1}]_{i=1,\ldots,\el;a_i \in \C^\times}$$ (see also \cite{DE}). We will
show that $\chi_q$ behaves well with respect to the restriction to
$\uqgg$ and to the quantum affine subalgebras of $\uqg$. We call
$\chi_q(V)$ the $q$--character of representation $V$. We hope that the
$q$--characters could be used more efficiently than their classical
counterparts in describing irreducible $\uqg$--modules and their
tensor products. In particular, we will define the {\em screening
operators} $S_i$ on $\yy$ and conjecture that the image of the
homomorphism $\chi_q$ equals the intersection of the kernels of $S_i,
i=1,\ldots,\el$ (we prove this for $\g=\sw_2$). Because of this
conjecture, we expect that the $q$--character of the irreducible
$\uqg$--module with a given highest weight can be found in a purely
combinatorial way.

The motivation for our construction of the screening operators comes
from our theory of deformed $\W$--algebras
\cite{FR:crit,FR:simple}. Let us first briefly describe the analogous
picture in the $q=1$ case. There are essentially three equivalent
definitions of the classical undeformed $\W$--algebra: as the result
of the Drinfeld-Sokolov reduction \cite{DS}, as the center of the
completed enveloping algebra $U(\G)$ at the critical level
\cite{FF:ds}, and as the algebra of integrals of motion of the Toda
field theory \cite{FF:laws}. According to the last definition, the
$\W$--algebra is a subalgebra in a Heisenberg algebra, which equals
the intersection of the kernels of the screening operators. This
definition can be thought of as an explicit description of the center
of $U(\G)$ at the critical level (see \cite{FF:ds}).

Now consider the $q$--version of this picture. We expect that the
center of $\uqg$ at the critical level is isomorphic to the
$q$--deformed classical $\W$--algebra and hence can be described as
the intersection of kernels of screening operators. But for $q\neq
1$, we have an injective homomorphism from $\on{Rep} \uqg$ to the
center of $\uqg((z))$ \cite{RS,DE}, and this allows us to view
$\on{Rep} \uqg$ as the ``space of fields'' of the $q$--deformed
$\W$--algebra. This suggests to us that $\on{Rep} \uqg$ can also be
described using the screening operators, which are precisely the
operators $S_i$.

In the course of writing this paper, we learned from M.~Kashiwara
about the work of H.~Knight \cite{Kn}. Knight proposed a character
theory for finite-dimensional representations of Yangians using a
different approach. Some further results in this direction were
subsequently obtained by Chari and Pressley \cite{CP5}. Knight's
results can be carried over to the case of quantum affine
algebras. However, he was unable to show the multiplicative property
of characters in this case because a certain result on the structure
of comultiplication in $\uqg$ was not available at the time. This
result is available now (see \lemref{comult}), and using it one can
extend Knight's construction to the case of quantum affine
algebras. In hindsight, it turns out that Knight's characters
essentially coincide with our $q$--characters. One of the advantages
of our definition is that the multiplicative property of the
$q$--characters mentioned above follows automatically, because
$\chi_q$ is manifestly a ring homomorphism.

The $q$--characters are also closely related to the formulas for the
spectra of transfer-matrices in integrable spin chains associated to
$\uqg$, obtained by the analytic Bethe Ansatz
\cite{Ba,Re1,Re2,BR,KS1}. Our definition of $q$--characters allows us
to streamline the rather {\em ad hoc} method of writing formulas for
these eigenvalues that has been used before (see \secref{Bethe}).

The results of this paper can be generalized in a straightforward way
to the Yangians and the twisted quantum affine algebras. Furthermore,
it turns out that $\on{Rep} \uqg$ and $\on{Rep} U_t(\GL)$, where $\GL$
is the Langlands dual affine Kac-Moody algebra to $\G$ (it is twisted
if $\g$ is non-simply laced) are two different classical limits of the
quantum deformed $\W$--algebra $\W_{q,t}(\g)$. Thus, the fields from
$\W_{q,t}(\g)$ may be viewed as $(q,t)$--characters, which are
simultaneous quantizations of the $q$--characters of $\uqg$ and the
$t$--characters of $U_t(\GL)$. The former appear when $t \arr 1$, and
the latter appear when $q \arr \exp (\pi i/r^\vee)$ (see
\cite{FR:simple}). The precise nature of this duality deserves further
study.

The paper is arranged as follows. In Sect.~2 we give the main
definitions and recall some of the results of Chari and Pressley on
finite-dimensional representations of $\uqg$. We define the
$q$--characters in Sect.~3. In Sects.~4 and 5 we present some results
and a conjecture on the structure of $q$--characters. In
particular, we discuss the connection between the $q$--characters and
the $R$--matrices. In Sect.~6 we explain the formulas for the spectra
of transfer-matrices obtained by the Bethe Ansatz from the point of
view of the $q$--characters. In Sect.~7 we define the screening
operators and state the conjecture characterizing the
$q$--characters. In Sect.~8 we motivate this conjecture by providing a
detailed analysis of the connection between $\on{Rep} \uqg$ and the
deformed $\W$--algebras.

\smallskip

\noindent{\bf Acknowledgments} We would like to thank J.~Beck,
V.~Chari, D.~Kazhdan, Ya.~Soi\-bel\-man, and especially M.~Kashiwara
for valuable comments and references. Since the time when the first
version of this paper was released, E.~Mukhin has found a
counter-example to the positivity conjecture stated there. We thank
him for his permission to include it in this version (although we
certainly wish he had not found it!), and for his careful reading of
the paper.

Some of the results of this paper have been reported by the first
author at the conferences at UC Riverside in December 1997, NCSU
Raleigh in May 1998, and ETH Z\"urich in June 1998.

The research of E.F. was supported by a grant from the Packard
Foundation and the research of N.R. was supported by an NSF grant.

\section{Background}

\subsection{Root data}    \label{cartan}

Let $\g$ be a simple Lie algebra of rank $\el$. We denote by $I$ the
set $\{ 1,\ldots,\el \}$. Let $h^\vee$ be the Coxeter number of
$\g$. Let $\langle \cdot,\cdot \rangle$ be the invariant inner product
on $\g$, normalized as in \cite{Kac}, so that the square of length of
the maximal root equals $2$ with respect to the induced inner product
on the dual space to the Cartan subalgebra $\h$ of $\g$ (also denoted
by $\langle \cdot,\cdot \rangle$). Let $\{ \al_1,\ldots,\al_\el \}$
and $\{ \om_1,\ldots,\om_\el \}$ be the sets of simple roots and of
fundamental weights of $\g$, respectively. We have:
$$
\langle \al_i,\om_j \rangle = \frac{\langle \al_i,\al_i \rangle}{2}
\delta_{i,j}.
$$
Let $r^\vee$ be the maximal number of edges connecting two vertices of
the Dynkin diagram of $\g$. Thus, $r^\vee=1$ for simply-laced $\g$,
$r^\vee=2$ for $B_\el, C_\el, F_4, G_2$, and $r^\vee=3$ for $D_4$.

From now on we will use the inner product
$$
(\cdot,\cdot) = r^\vee \langle \cdot,\cdot \rangle
$$
on $\h^*$. Set
$$
D = \on{diag}(\rr_1,\ldots,\rr_\el),
$$
where
\begin{equation}    \label{di}
\rr_i = \frac{(\al_i,\al_i)}{2} = r^\vee \frac{\langle \al_i,\al_i
\rangle}{2}.
\end{equation}
All $\rr_i$'s are integers; for simply-laced $\g$, $D$ is the identity
matrix.

Now let $C = (C_{ij})_{1\leq i,j\leq \el}$ and $(I_{ij})_{1\leq
i,j\leq \el}$ be the {\em Cartan matrix} and the {\em incidence matrix}
of $\g$, respectively, so that $C = 2 \on{Id} - I$. We have:
$$
C_{ij} = \frac{2(\al_i,\al_j)}{(\al_i,\al_i)}.
$$
Let $B = (B_{ij})_{1\leq i,j\leq \el}$ be the symmetric matrix
$$
B = D C,
$$
i.e.,
$$
B_{ij} = (\al_i,\al_j) = r^\vee \langle \al_i,\al_j \rangle.
$$

Let $q \in \C^\times$ be such that $|q|<1$. Set $q_i = q^{r_i}$.

We will use the standard notation
$$
[n]_q = \frac{q^n - q^{-n}}{q - q^{-1}}.
$$
Introduce the $\el \times \el$ matrices $B(q), C(q), D(q)$ by the
formulas
\begin{align*}
B_{ij}(q) &= [B_{ij}]_q,\\
C_{ij}(q) &= (q^{r_i} + q^{-r_i}) \delta_{ij} - [I_{ij}]_q,\\
D_{ij}(q) &= \delta_{ij} [r_i]_q.
\end{align*}
We have:
$$
B(q) = D(q) C(q).
$$

\subsection{Quantum affine algebras}

\begin{definition}[\cite{Dr1,J}] Let $\uqg$ be the
associative algebra over $\C$ with generators $x_i^{{}\pm{}}$,
$k_i^{{}\pm 1}$ ($i=0,\ldots,\el$), and relations:
\begin{align*}
k_ik_i^{-1} = k_i^{-1}k_i &=1,\quad \quad k_ik_j =k_jk_i,\\
k_ix_j^{{}\pm{}}k_i^{-1} &= q^{{}\pm B_{ij}}x_j^{{}\pm},\\ [x_i^+ ,
x_j^-] &= \delta_{ij}\frac{k_i - k_i^{-1}}{q_i -q_i^{-1}},\\
\sum_{r=0}^{1-C_{ij}}\left[\begin{array}{cc} 1-C_{ij} \\ r \end{array}
\right]_{q_i}
(x_i^{{}\pm{}})^rx_j^{{}\pm{}}&(x_i^{{}\pm{}})^{1-C_{ij}-r} =0, \ \ \
\ i\ne j.
\end{align*}

We introduce a $\Z$-gradation on $\uqg$ by the formulas: $\deg x_0^\pm
= \pm 1, \deg x_i^\pm = \deg k_i = 0, i \in I=\{ 1,\ldots,\el\}$.

Denote the subalgebra of $\uqg$ generated by $k_i^{\pm 1}, x_i^+$
(resp., $k_i^{\pm 1}, x_i^-$), $i=0,\ldots\el$, by $U_q \bb_+$ (resp.,
$U_q \bb_-$).

The algebra $\uqgg$ is defined as the subalgebra of $\uqg$ with
generators $x_i^{{}\pm{}}$, $k_i^{{}\pm 1}$, where $i\in I$.

$\uqg$ has a structure of a Hopf algebra with the comultiplication
given on generators by:
\begin{align*}
\Delta(k_i) &= k_i \otimes k_i, \\
\Delta(x^+_i) &= x^+_i \otimes 1 + k_i \otimes x^+_i,\\
\Delta(x^-_i) &= x^-_i \otimes k_i^{-1} + 1 \otimes x^-_i,
\end{align*}
\end{definition}

\begin{rem}    \label{choice}
This comultiplication differs from the original one \cite{Dr1,J}. The
difference is essentially accounted for by the automorphism of $\uqg$
sending $k_i$ to $k_i^{-1}$, $x^\pm_i$ to itself and $q$ to
$q^{-1}$. The reason for this choice of comultiplication is that it is
in terms of this comultiplication that the universal $R$--matrix has
the form given in \secref{rmatr}.\qed
\end{rem}

The following theorem describes the Drinfeld ``new'' realization of
$\uqg$ \cite{Dr}.

\begin{thm}[\cite{Dr,KhT,LSS,Beck}]
{\em The algebra $\uqg$ has another realization as the algebra with
generators $x_{i,n}^{\pm}$ ($i\in I=\{ 1,\ldots,\el\}$, $n\in\Z$),
$k_i^{\pm 1}$ ($i\in I=\{ 1,\ldots,\el\}$), $h_{i,n}$ ($i\in I$, $n\in
\Z\backslash\{0\}$) and central elements $c^{\pm 1/2}$, with the
following relations:
\begin{align*}
k_ik_j = k_jk_i,\quad & k_ih_{j,n} =h_{j,n}k_i,\\
k_ix^\pm_{j,n}k_i^{-1} &= q^{\pm B_{ij}}x_{j,n}^{\pm},\\
[h_{i,n} , x_{j,m}^{\pm}] &= \pm \frac{1}{n} [n B_{ij}]_q c^{\mp
{|n|/2}}x_{j,n+m}^{\pm},\\ x_{i,n+1}^{\pm}x_{j,m}^{\pm}
-q^{\pm B_{ij}}x_{j,m}^{\pm}x_{i,n+1}^{\pm} &=q^{\pm
B_{ij}}x_{i,n}^{\pm}x_{j,m+1}^{\pm}
-x_{j,m+1}^{\pm}x_{i,n}^{\pm},\\ [h_{i,n},h_{j,m}]
&=\delta_{n,-m} \frac{1}{n} [n B_{ij}]_q \frac{c^n -
c^{-n}}{q-q^{-1}},\\ [x_{i,n}^+ , x_{j,m}^-]=\delta_{ij} & \frac{
c^{(n-m)/2}\phi_{i,n+m}^+ - c^{-(n-m)/2} \phi_{i,n+m}^-}{q_i -
q_i^{-1}},\\
\sum_{\pi\in\Sigma_s}\sum_{k=0}^s(-1)^k\left[\begin{array}{cc} s \\
k \end{array} \right]_{q_i} x_{i, n_{\pi(1)}}^{\pm}\ldots
x_{i,n_{\pi(k)}}^{\pm} & x_{j,m}^{\pm} x_{i,
n_{\pi(k+1)}}^{\pm}\ldots x_{i,n_{\pi(s)}}^{\pm} =0,\ \ i\ne
j,
\end{align*}
for all sequences of integers $n_1,\ldots,n_s$, where $s=1-C_{ij}$,
$\Sigma_s$ is the symmetric group on $s$ letters, and
$\phi_{i,n}^{\pm}$'s are determined by the formal power series}
\begin{equation}    \label{series}
\sum_{n=0}^{\infty}\phi_{i,\pm n}^{\pm}u^{\pm n} = k_i^{\pm
1} \exp\left(\pm(q-q^{-1})\sum_{m=1}^{\infty}h_{i,\pm m} u^{\pm
m}\right).
\end{equation}
\end{thm}

\begin{rem} Note that the generators $h_{i,n}$ correspond to
$\dfrac{q_i-q_i^{-1}}{q-q^{-1}} {\mathcal H}_{i,n}$ of \cite{CP}.\qed
\end{rem}

Let $Q$ be the root lattice of $\G$. We introduce the $Q$--gradation
on $\uqg$ by the formulas: $\deg x^\pm_{i,n} = \pm \al_i, \deg k_i =
\deg h_{i,n} = \deg c^{1/2} = 0$.

For any $a \in \C^\times$, there is a Hopf algebra automorphism
$\tau_a$ of $\uqg$ defined on the generators by the following
formulas:
$$
\tau_a(x_{i,r}^{\pm})=a^rx_{i,r}^{\pm}, \quad \quad
\tau_a(\phi_{i,r}^{\pm})=a^r\phi_{i,r}^{\pm},
$$
$$
\tau_a(c^{1/2})=c^{1/2}, \quad \quad \tau_a(k_i)=k_i,
$$
for all $i=1,\ldots,\ell, r \in \Z$. Given a $\uqg$--module $V$ and $a
\in C^\times$, we denote by $V(a)$ the pull-back of $V$ under $\tau_a$.

Introduce new variables $\wt{k}_i^{\pm 1}, i\in I$, such that $$k_j =
\prod_{i\in I} \wt{k}_i^{C_{ij}}.$$ Thus, while $k_i$ corresponds to
the simple root $\al_i$, $\wt{k}_i$ corresponds to the fundamental
weight $\omega_i$. We extend the algebra $\uqg$ by replacing the
generators $k_i^{\pm 1}, i \in I$ with $\wt{k}_i^{\pm 1}, i\in
I$. From now on $\uqg$ will stand for the extended algebra.

\subsection{Finite-dimensional representations of $\uqg$}

In this section we recall some of the results of Chari and Pressley
\cite{CP,CP2,CP3,CP4}.

Let $P$ be the weight lattice of $\g$.  A representation $W$ of
$\uqgg$ is said to be of type 1 if it is the direct sum of its weight
spaces
$$W_{\lambda} = \{w\in W | k_i \cdot w = q^{(\lambda,\al_i)}w\}, \quad
\quad \lambda\in P.$$ If $W_\lambda\ne 0$, then $\lambda$ is called a
weight of $W$. A vector $w\in W_{\lambda}$ is called a highest weight
vector if $x_i^+ \cdot w =0$ for all $i\in I$, and $W$ is a highest
weight representation with highest weight $\lambda$ if $W=\uqgg \cdot
w$ for some highest weight vector $w\in W_\lambda$. In that case, $\la
\in P^+$, the set of dominant weights.

A representation $V$ of $\uqg$ is called of type 1 if $c^{1/2}$ acts
as the identity on $V$, and if $V$ is of type 1 as a representation of
$\uqgg$. A vector $v\in V$ is called a highest weight vector if
\begin{equation}    \label{hwv}
x_{i,r}^+ \cdot v=0,\quad \quad \phi_{i,r}^\pm \cdot v =
\psi_{i,r}^\pm v,\quad \quad c^{1/2} v =v,
\end{equation}
for some complex numbers $\psi_{i,r}^{\pm}$. A type 1
representation $V$ is a highest weight representation if $V=\uqg \cdot
v$, for some highest weight vector $v$. In that case the set
$(\psi_{i,r}^{\pm})_{i\in I,r\in\Z}$ is called the highest weight
of $V$.

If $\lambda\in P^+$, let ${\Cal P}^\lambda$ be the set of all
$I$--tuples $(P_i)_{i\in I}$ of polynomials $P_i\in\C[u]$, with
constant term 1, such that $\deg(P_i)=\lambda(\al_i^\vee)$ for all
$i\in I$. Such an $I$--tuple is then said to have degree $\la$. Set
${\Cal P}=\bigcup_{\lambda\in P^+}{\Cal P}^\lambda$.

Set
$$
\Phi_i^\pm(u) = \sum_{n=0}^\infty \phi_{i,n}^\pm u^{\pm n}, \quad
\quad \Psi_i^\pm(u) = \sum_{n=0}^\infty \psi_{i,n}^\pm u^{\pm n}.
$$

\begin{thm}[\cite{CP,CP3}]    \label{ChP2}
\hfill

{\em (1) Every finite-dimensional irreducible
representation of $\uqg$ can be obtained from a type 1 representation
by twisting with an automorphism of $\uqg$.

(2) Every finite-dimensional irreducible representation of $\uqg$ of
type 1 is a highest weight representation.

(3) Let $V$ be a finite-dimensional irreducible representation of
$\uqg$ of type 1 and highest weight $(\Psi_{i,r}^{\pm})_{i\in
I,r\in\Z}$. Then, there exists $\bp=(P_i)_{i\in I}\in\calp$ such
that
\begin{equation}
\Psi_i^\pm(u) = q_i^{\deg(P_i)}\frac{P_i(uq_i^{-1})}{P_i(uq_i)}.
\end{equation}
as an element of $\C[[u^{\pm 1}]]$.

Assigning to $V$ the set $\bp$ defines a bijection between $\calp$ and
the set of isomorphism classes of finite-dimensional irreducible
representations of $\uqg$ of type 1. The irreducible representation
associated to $\bp$ will be denoted by $V(\bp)$.

(4) If $\bp=(P_i)_{i\in I}\in\calp$, $a\in\C^\times$, and if
$\tau_a^*(V(\bp))$ denotes the pull-back of $V(\bp)$ by the automorphism
$\tau_a$, we have
$$\tau_a^*(V(\bp))\cong V(\bp^a)$$
as representations of $\uqg$, where $\bp^a=(P_i^a)_{i\in I}$ and
$$P_i^a(u)=P_i(ua).$$

(5) Let $\bp$, $\bq\in\calp$ be as above, and let $v_{\bp}$ and
$v_{\bq}$ be highest weight vectors of $V(\bp)$ and $V(\bq)$,
respectively. Denote by $\bp\ot\bq$ the $I$-tuple $(P_iQ_i)_{i\in
I}$. Then $V(\bp\ot\bq)$ is isomorphic to a quotient of the
subrepresentation of $V(\bp)\ot V(\bq)$ generated by the tensor
product of the highest weight vectors.}
\end{thm}

\begin{rem} An analogous classification result for Yangians has been
obtained earlier by Drinfeld \cite{Dr}. Because of that, the
polynomials $P_i(u)$ are sometimes called Drinfeld polynomials.\qed
\end{rem}

\begin{rem} Note that in our notation the polynomial $P_i(u)$
corresponds to the polynomial $P_i(uq_i^{-1})$ in the notation of
\cite{CP,CP3}.\qed
\end{rem}

For each $i\in I$ and $a \in \C^\times$, we define the irreducible
representation $V_{\omega_i}(a)$ as $V(\bp^{(i)}_a)$, where
$\bp^{(i)}_a$ is the $I$--tuple of polynomials, such that
$P_i(u)=1-ua$ and $P_j(u)=1, \forall j \neq i$. We call
$V_{\omega_i}(a)$ the $i$th fundamental representation of $\uqg$. Note
that in general $V_{\omega_i}(a)$ is reducible as a $\uqgg$--module.

\thmref{ChP2} implies the following

\begin{cor}[\cite{CP3}]    \label{generated}
Any irreducible finite-dimensional representation of $\uqg$ occurs as
a quotient of the submodule of the tensor product
$V_{\omega_{i_1}}(a_1) \otimes \ldots \otimes V_{\omega_{i_n}}(a_n)$,
generated by the tensor product of the highest weight vectors. The
parameters $(\omega_{i_1},a_1),\ldots,(\omega_{i_n},a_n)$ are uniquely
determined by this representation up to permutation.
\end{cor}

\begin{rem} For $U_q \sun$, V.~Ginzburg and E.~Vasserot have given an
alternative geometric construction of irreducible finite-dimensional
representations \cite{GV,V}.\qed
\end{rem}

\subsection{Spectra of $\Phi^\pm(u)$ on finite-dimensional
representations}

It follows from the defining relations that the operators
$\phi_{i,n}^\pm$ commute with each other. Hence we can decompose any
representation $V$ of $\uqg$ into a direct sum $V = \oplus
V_{(\ga^\pm_{i,n})}$, where
$$
V_{(\ga^\pm_{i,n})} = \{ x \in V | (\phi^\pm_{i,n} - \ga^\pm_{i,n})^p
\cdot x = 0, \on{for} \on{some} p, \forall i,n \}.
$$
Given a collection $(\ga^\pm_{i,n})$ of eigenvalues, we form the
generating functions
$$
\ga^\pm_i(u) = \sum_{n>0} \ga^\pm_{i,\pm n} u^{\pm n}.
$$
We will refer to a series $\ga^\pm_i(u)$ occurring on a given
representation $V$ as an eigenvalue of $\Phi^\pm_i(u)$ on $V$. The
following result is a generalization of \thmref{ChP2}.

\begin{prop}    \label{spectra}
The eigenvalues $\ga^\pm_i(u)$ of $\Phi_i^\pm(u)$ on any
finite-dimensional representation of $\uqg$ have the form:
\begin{equation}    \label{form}
\ga_i^\pm(u) = q_i^{\deg Q_i - \deg R_i} \frac{Q_i(uq_i^{-1})
R_i(uq_i)}{Q_i(uq_i) R_i(uq_i^{-1})},
\end{equation}
as elements of $\C[[u]]$ and $\C[[u^{-1}]]$, respectively, where
$Q_i(u), R_i(u)$ are polynomials in $u$ with constant term $1$.
\end{prop}

\begin{proof}
Let $\uqg_{\{i\}}$ be the subalgebra of $\uqg$ generated by $k_i^{\pm 1},
h_{i,n}, x^\pm_{i,n}, n \in \Z$. This subalgebra is isomorphic to
$U_{q_i} \widehat{\sw}_2$. The eigenvalues of $\Phi^+_i(u)$ on a
finite-dimensional representation $V$ of $\uqg$ coincide with the
eigenvalues of $\Phi^\pm_i(u)$ on the restriction of $V$ to
$\uqg_{\{i\}}$. Hence it suffices to prove the statement of the proposition
when $\g=\sw_2$.

The finite-dimensional irreducible representations of $\uqsl$ have
been classified in \cite{CP2} (see also \cite{CP3}). The result is
recalled in \thmref{classl2} below. According to this result, each
irreducible representation of $\uqsl$ is isomorphic to the tensor
product $W_{r_1}(b_1) \otimes \ldots \otimes W_{r_m}(b_m)$. The
representation $W_i(b)$ is defined in \secref{casesl2} and the
eigenvalues of $\Phi^\pm(u)$ on it are given by formula
\eqref{spec2}. These eigenvalues are in the form \eqref{form}. Hence
it suffices to show that the eigenvalues of $\Phi^\pm(u)$ on the
tensor product $V \otimes W$ are equal to the products of the
eigenvalues on $V$ and $W$. This follows from the formula for the
comultiplication of $h_{i,n}$, which is proved in \cite{Da}. Actually,
in \cite{Da} a more precise formula is proved; that formula has also
appeared in \cite{JKK,BCP}. Note that we need this formula only in the
case $\g=\sw_2$.

\begin{lem}    \label{comult}
On representations of $\uqg$ on which $c$ acts as the identity,
\begin{equation}    \label{com}
\Delta(h_{i,\pm n}) = h_{i,\pm n} \otimes 1 + 1 \otimes h_{i,\pm n} +
\on{terms} \; \on{in} \; U_\mp \otimes U_\pm,\ \quad \quad n>0,
\end{equation}
where $U_+$ (resp., $U_-$) is the subalgebra of $\uqg$ spanned by
elements of positive (resp., negative) $Q$--degree.
\end{lem}

Let us concentrate on the case of $\Phi^+(u)$ and hence $h_n,
n>0$. The case of $\Phi^-(u)$ is analyzed in the same way. Let $V$ and
$W$ be two representations of $\uqsl$. We can decompose $V$ and $W$
into direct sums $V = \oplus_{p \in \Z} V_p, W = \oplus_{p \in \Z}
W_p$, where
$$
V_p = \{ x \in V | k \cdot x = q^p x \}.
$$
Since $k$ commutes with all $h_n$, all $V_p, W_p$ are $h_n$--invariant
for all $n \in \Z$. We can choose bases $\{ v_{\al}^p \}$ of $V_p$ and
$\{ w_{\beta}^p \}$ in $W$ with respect to which $h_n$'s are
upper-triangular for all $n$. Let us now order the basis $\{ v_{\al}^p
\otimes w_{\beta}^q \}$ of $V \otimes W$ in such a way that $v_\al^p
\otimes w_\beta^q < v_{\al'}^r \otimes w_{\beta'}^s$ if either $q>s$
or $q=s$ and $p<r$. For fixed $p,q,r,s$ we keep the old orderings on
$\al,\beta,\al',\beta'$. Formula \eqref{com} shows that $\Delta(h_n)$
is upper triangular in this basis of $V \otimes W$, and the
eigenvalues of $\Delta(h_n)$ are equal to the sums of the eigenvalues
of $h_n$ on $V$ and $W$. Using formula \eqref{series} for
$\Phi^+(u)$ and the formula $\Delta(k_i) = k_i \otimes k_i$, we
obtain that the eigenvalues of $\Phi^+(u)$ on $V \otimes W$ are
products of its eigenvalues on $V$ and $W$, which is what we needed to
show.
\end{proof}

\begin{rem} The statement analogous to \propref{spectra} in the case
of the Yangians has been proved by Knight \cite{Kn}, Prop.~4. Our
proof is similar to his proof.

In the same way as above, one can derive from \lemref{comult} that the
eigenvalues of $\Phi^\pm_i(u)$ on the tensor product $V \otimes W$ are
the products of its eigenvalues on $V$ and $W$ for any $\uqg$.\qed
\end{rem}

\section{Definition of $q$--characters}

\subsection{Transfer-matrices}    \label{tm}

The completed tensor product $\uqg \;\widehat{\otimes}\; \uqg$ contains a
special element $\R$ called the universal $R$--matrix (at level
$0$). It actually lies in $U_q \bb_+ \;\widehat{\otimes}\; U_q \bb_-$ and
satisfies the following identities:
\begin{align}    \label{r1}
\Delta'(x) &= \R \Delta(x) \R^{-1}, \quad \quad \forall x
\in \uqg,\\    \label{r2}
(\Delta \otimes \on{id}) \R &= \R^{13} \R^{23}, \quad
\quad (\on{id} \otimes \Delta) \R = \R^{13} \R^{12}.
\end{align}
Note that the last two equations imply that $\R$ satisfies the
Yang-Baxter equation. For more details, see \cite{Dr2,EFK}.

Now let $(V,\pi_V)$ be a finite-dimensional representation of
$\uqg$. Define the $L$--{\em operator} corresponding to $V$ by the
formula
\begin{equation}    \label{lv}
L_V = L_V(z) = (\pi_{V(z)} \otimes \on{id})(\R).
\end{equation}

Now define the {\em transfer-matrix} corresponding to $V$ by
\begin{equation}    \label{tv}
t_V = t_V(z) = \on{Tr}_V \; q^{2\rho} \; L_V(z),
\end{equation}
where by definition $q^{2\rho} = \wt{k}_1^2 \ldots \wt{k}_\el^2$.
Then for each $V$,
$$
t_V(z) = \sum_{n\leq 0} t_V[n] z^{-n},
$$
where $t_V[n] \in U_q \bb_-$. Furthermore, $t_V[n]$ has degree $n$
with respect to the $\Z$--gradation on $U_q \bb_-$ (see Definition
1). Thus, for each $x \in \C^\times$, $t_V(x)$ lies in the completion
$\wt{U}_q \bb_-$ of $U_q \bb_-$ corresponding to this gradation.

\begin{lem}
\hfill

(1) For any pair of finite-dimensional representations, $V$ and $W$,
$$[t_V(z),t_W(w)] = 0;$$

(2) Given a short exact sequence $0 \arr V \arr W \arr U \arr 0$,
$t_W(z) = t_V(z) + t_U(z)$;

(3) $t_{V \otimes W}(z) = t_V(z)t_W(z)$;

(4) $t_{V(a)}(z) = t_V(za), \forall a \in \C^\times$.
\end{lem}

\begin{proof}
Parts (2)--(4) follow immediately from the definition of $t_V(z)$.

It is clear that $$R_{V,W}(z,w) = (\pi_{V(z)} \otimes
\pi_{W(w)})(\R)$$ is a well-defined element of $\on{End}(V \otimes
W)[[z,w^{-1}]]$. The Yang-Baxter equation gives:
$$
L_V(z) L_W(w) = R_{V,W}(z,w)^{-1} L_W(w) L_V(z) R_{V,W}(z,w).
$$
Taking the traces and using (2) we obtain: $t_V(z) t_W(w) = t_W(w)
t_V(z)$ as formal power series.
\end{proof}

The Lemma implies the following statement.

\begin{prop}    \label{nuq}
The linear map $\nu_q$ sending $V \in \on{Rep} \uqg$ to $t_V(z) \in
\uqb[[z]]$ is a $\C^\times$--equivariant ring homomorphism from
$\on{Rep} \uqg$ onto a commutative subalgebra of $\uqb[[z]]$.
\end{prop}

\subsection{Analogue of the Harish-Chandra homomorphism}    \label{hc}

Denote by $\uqbb$ the subalgebra of $\uqg$ generated by $x^\pm_{i,n},
k_i, h_{i,n}, i \in I, n\leq 0$. As a vector space, $\uqbb$ can be
decomposed as follows: $\uqbb = \uqn_- \otimes \uqh \otimes \uqn_+$,
where $\uqn_\pm$ (resp., $\uqh$) is generated by $x^\pm_{i,n}$ (resp.,
$k_i, h_{i,n}, i \in I, n\leq 0$). Hence
\begin{equation}    \label{deco}
\uqbb = \uqh \oplus (\uqbb \cdot (\uqn_+)_0) \oplus ((\uqn_-)_0 \cdot
\uqbb),
\end{equation}
where $(\uqn_\pm)_0$ stands for the augmentation ideal of
$\uqn_\pm$. Denote by ${\mathbf h}_q$ the projection $\uqbb \arr \uqh$
along the last two summands. We denote by the same letter its
restriction to $\uqb$.

\begin{definition}
The map $\chi_q: \on{Rep} \uqg \arr \uqh[[z]]$ is the composition of
$\nu_q: \on{Rep} \uqg \arr \uqb[[z]]$ and ${\mathbf h}_q: \uqb \arr \uqh$.
\end{definition}

\begin{lem}
The map $\chi_q$ is a ring homomorphisms.
\end{lem}

\begin{proof}
Let ${\mathfrak z}_q(\G)$ be the commutative subalgebra of $\uqbb$
generated by $t_V[n], V \in \on{Rep} \uqg$, $n\leq 0$. Let us show
that the restriction of ${\mathbf h}_q$ to ${\mathfrak z}_q(\G)$ is a
ring homomorphism ${\mathfrak z}_q(\G) \arr \uqh$. Since $\nu_q$ is a
ring homomorphism according to \propref{nuq}, this will prove the
statement of the lemma.

Consider two elements $A, B \in {\mathfrak z}_q(\G)$, By construction,
both of them have degree $0$ with respect to the $Q$--gradation on
$\uqg$. Hence we can write $A = A_0 + A_1$, where $A_0 = {\mathbf
h}_q(A) \in \uqh$, $A_1 \in (\uqn_-)_0 \cdot \uqbb$, and $B = B_0 +
B_1$, where $B_0 = {\mathbf h}_q(B) \in \uqh$, $B_1 \in \uqbb \cdot
(\uqn_+)_0$. But then ${\mathbf h}_q(AB) = {\mathbf h}_q(A_0 B_0) +
{\mathbf h}_q(A_0 B_1) + {\mathbf h}_q(A_1B_0) + {\mathbf h}_q(A_1
B_1) = A_0 B_0$, which is what we needed to prove.
\end{proof}

\subsection{The image of $\chi_q$}    \label{rmatr}

In order to describe the image of $\chi_q$, we need an explicit
formula for the universal $R$--matrix. This formula was obtained by
Khoroshkin and Tolstoy \cite{KhT} and, by a different method, by
Levendorsky-Soibelman-Stukopin \cite{LSS} (in the case of $\uqsl$) and
Damiani \cite{Da} (for general $\uqg$).

Denote by $\uqnn_\pm$ the subalgebra of $\uqg$ generated by
$x^\pm_{i,n}, i\in I, n\in \Z$. Let $\wt{B}(q)$ be the inverse matrix
to $B(q)$ from \secref{cartan}. The formula for the universal
$R$--matrix then reads:
\begin{equation}    \label{ktr}
\R  = \R^+ \R^0 \R^- T,
\end{equation}
where
\begin{equation}
\R^0 = \exp \left( - (q-q^{-1}) \sum_{n>0} \frac{n}{[n]_q}
\wt{B}_{ij}(q^n) h_{i,n} \otimes h_{j,-n} \right),
\end{equation}
$\R^\pm \in \uqnn_\pm \otimes \uqnn_\mp$, and $T$ acts as follows: if
$x,y$ satisfy $k_i \cdot x = q^{(\la,\al_i)} x, k_i \cdot y =
q^{(\mu,\al_i)} y$, then
\begin{equation}    \label{t}
T \cdot x \otimes y = q^{-(\la,\mu)} x \otimes y.
\end{equation}

\begin{rem} Note that the above formula for $T$ differs from Drinfeld's
formula \cite{Dr2} by the replacement $q \arr q^{-1}$. This is because
the comultiplication that we use differs in the same way from that of
\cite{Dr2} (cf. \remref{choice}).\qed
\end{rem}

The computation of $\chi_q(V)$ for $V \in \on{Rep} \uqg$ proceeds
along the lines of the computation of Ding and Etingof in Sect.~3 of
\cite{DE}. First, the projection onto $\uqh$ eliminates the factor
$\R^-$ from the formula, and then taking the trace eliminates $\R^+$
(recall that $\uqn_+$ acts nilpotently on $V$). Hence we are left with
\begin{multline}    \label{muqv}
\chi_q(V) = \\ \on{Tr}_V \; \left[ q^{2\rho} \exp \left( - (q-q^{-1})
\sum_{n>0} \frac{n}{[n]_q} \wt{B}_{ij}(q^n) z^n \pi_V(h_{i,n}) \otimes
h_{j,-n} \right) (\pi_V \otimes 1)(T) \right].
\end{multline}

Now let
\begin{equation}    \label{connection}
\wt{h}_{i,-m} = \sum_{j\in I} \frac{q_i^m - q_i^{-m}}{q^m -
q^{-m}} \wt{B}_{ij}(q^m)
h_{j,-m} = \sum_{j\in I} \wt{C}_{ji}(q^m) h_{j,-m},
\end{equation}
where $\wt{C}(q)$ is the inverse matrix to $C(q)$ defined in
\secref{cartan}. Set
\begin{equation}    \label{Yia}
Y_{i,a} = q^{2(\rho,\omega_i)} \wt{k}_i^{-1} \exp \left( - (q-q^{-1})
\sum_{n>0} \wt{h}_{i,-n} z^n a^n \right), \quad \quad a \in \C^\times.
\end{equation}
We assign to $Y_{i,a}^{\pm 1}$ the weight $\pm \omega_i$.

Before we state our theorem, we need to introduce some more notation.

Let $$\chi: \on{Rep} \uqgg \arr \Z[y_i^{\pm 1}]_{i\in I}$$ be the
ordinary character homomorphism, $\beta$ be the homomorphism
$$\Z[Y_{i,a_i}^{\pm 1}]_{i\in I; a_i \in \C^\times} \arr \Z[y_i^{\pm
1}]_{i\in I}$$ sending $Y_{i,a_i}^{\pm 1}$ to $y_i^{\pm 1}$, and
$$
\on{res}: \on{Rep} \uqg \arr \on{Rep} \uqgg
$$
be the restriction homomorphism.

Given a subset $J$ of $I$ we denote by $\uqg_J$ the subalgebra of
$\uqg$ generated by $k_i^{\pm 1}, h_{i,n}, x^\pm_{i,n}, i\in J, n \in
\Z$. Let $$\on{res}_J: \on{Rep} \uqg \arr \on{Rep} \uqg_J$$ be the
restriction map and $\beta_J$ be the homomorphism $\Z[Y_{i,a_i}^{\pm
1}]_{i\in I} \arr \Z[Y_{i,a_i}^{\pm 1}]_{i\in J}$, sending
$Y_{i,a}^{\pm 1}$ to itself for $i \in J$ and to $1$ for $i
\not{\hspace*{-1mm}\in} J$.

\begin{thm}    \label{mainthm}
\hfill

{\em (1) $\chi_q$ is an injective homomorphism from $\on{Rep} \uqg$ to

$\Z[Y_{i,a_i}^{\pm 1}]_{i\in I; a_i \in \C^\times} \subset \uqh[[z]]$.

(2) The diagram
$$\begin{CD}
\on{Rep} \uqg  @>{\chi_q}>>  \Z[Y_{i,a_i}^{\pm 1}]_{i\in I}\\
@VV{\on{res}}V     @VV{\beta}V\\
\on{Rep} \uqgg @>{\chi}>>  \Z[y_i^{\pm 1}]_{i\in I}
\end{CD}$$

is commutative

(3) The diagram
$$\begin{CD}
\on{Rep} \uqg @>{\chi_q}>>  \Z[Y_{i,a_i}^{\pm 1}]_{i\in I}\\
@VV{\on{res}_J}V     @VV{\beta_J}V\\
\on{Rep} \uqg_J  @>{\chi_{q,J}}>>  \Z[Y_{i,a_i}^{\pm 1}]_{i\in J}
\end{CD}$$
is commutative.}
\end{thm}

\begin{proof} According to \propref{spectra}, the eigenvalues of
$$
\pi_V(k_i) \exp \left( (q-q^{-1}) \sum_{n>0} z^n \pi_V(h_{i,n})
\right) = \pi_V(\Phi^+_i(z))
$$
are given by
$$
q_i^{\deg Q_i - \deg R_i} \frac{Q_i(zq_i^{-1})
R_i(zq_i)}{Q_i(zq_i) R_i(zq_i^{-1})},
$$
where
$$
Q_i(z) = \prod_{r=1}^{k_i} (1-za_{ir}), \quad \quad R_i(z) =
\prod_{s=1}^{l_i} (1-zb_{is}).
$$
This implies that a typical eigenvalue of $\pi_V(h_{i,n})$ equals
$$\frac{q_i^n-q_i^{-n}}{n(q-q^{-1})} \left( \sum_{r=1}^{k_i} a_{ir}^n
- \sum_{s=1}^{l_i} b_{is}^n \right), \quad \quad n>0.$$

Substituting this into formula \eqref{muqv} we obtain that $\chi_q(V)$
is a linear combination of monomials
$$
\prod_{i\in I} \prod_{r=1}^{k_i} Y_{i,a_{ir}} \prod_{s=1}^{l_i}
Y_{i,b_{is}}^{-1},
$$
with positive integral coefficients. Furthermore, it follows from the
construction that the set of weights of these monomials is the set of
weights of $V$ counted with multiplicities. Thus, the image of
$\chi_q$ lies in $\Z[Y_{i,a_i}^{\pm 1}]_{i\in I; a_i \in \C^\times}$,
and we obtain part (2) of the theorem. Part (3) is also clear.

It remains to show that $\chi_q$ is injective. Note that
$\chi_q(V_{\omega_i}(a_i))$ equals $Y_{i,a_i}$ plus the sum of
monomials of lower weight. By \corref{generated}, for any
finite-dimensional irreducible representation $V$, $\chi_q(V)$ equals
$Y_{i_1,a_1} \ldots Y_{i_n,a_n}$, where the set
$(i_1,a_1),\ldots,$ $(i_n,a_n)$ is uniquely determined by $V$ up to
permutation, plus the sum of monomials of lower weight. Since
$Y_{i,a_i}$'s are algebraically independent in $\uqh[[z]]$, this shows
that $\chi_q$ is injective.
\end{proof}

\begin{cor}    \label{commring}
$\on{Rep} \uqg$ is a commutative ring that is isomorphic to
$\Z[t_{i,a_i}]_{i\in I,a_i \in \C^\times}$, where $t_{i,a}$ is the
class of $V_{\omega_i}(a)$.
\end{cor}

\begin{rem} D.~Kazhdan and V.~Chari have communicated to us two
alternative proofs of commutativity of $\on{Rep} \uqg$.

The first is based on the fact that for any pair of finite-dimensional
representations, $V$ and $W$, $PR_{V,W}(z): V(z) \otimes W \arr W
\otimes V(z)$, where $R_{V,W}(z) = (\pi_{V(z)} \otimes \pi_W)(\R) \in
\on{End}(V \otimes W)[[z]]$ and $P(a \otimes b) = b \otimes a$, is an
expansion of a meromorphic function in $z$, which is an isomorphism
for generic $z \in \C$.

The second proof relies on Proposition 5.1(b) from \cite{CP4}, which
states that $V(\bp)^*$ is isomorphic to $V(\bp^*)$, where $P^*_i(u) =
P_{\ovl{i}}(uq^\kappa)$, and $\kappa$ is a constant (actually,
$\kappa=-r^\vee h^\vee$). Here $\ovl{i}$ is defined by $\al_{\ovl{i}} =
- w_0(\al_i)$, where $w_0$ is the longest element of the Weyl group of
$\g$. Since $(V \otimes W)^* = W^* \otimes V^*$, one can use this
result to compare the composition factors in tensor products $V
\otimes W$ and $W \otimes V$.\qed
\end{rem}

\begin{rem} The reader may wonder why we do not consider the ``naive''
character $\on{Tr}_V \Phi^\pm_i(u)$. In the case of $\uqsl$, explicit
computation shows that $\on{Tr}_V \Phi^\pm(u)$ is independent of $u$,
and equals $\on{Tr}_V k^{\pm 1}$. This implies that in general
$\on{Tr}_V \Phi^\pm_i(u) = \on{Tr}_V k_i^{\pm 1}$.\qed
\end{rem}

\begin{rem} The definition of $q$--characters carries over to the case
of Yangians in a straightforward way. The resulting characters
essentially coincide with the characters introduced by Knight
\cite{Kn}.\qed
\end{rem}

\section{The structure of $q$--characters}

\subsection{The case of $\uqsl$}    \label{casesl2}

Let us recall the classification of finite-dimensional representations
of $\uqsl$ due to Chari and Pressley \cite{CP2,CP3}.

For each $r \in \Z, r>0$, and $a\in \C^\times$ set
$$
P(u)^{(r)} = \prod_{k=1}^r (1-uaq^{r-2k+1}).
$$
Denote by $W_r(a)$ the irreducible representation of $\uqsl$ with
highest weight $P(u)^{(r)}$. Recall that such a representation is
unique up to isomorphism. These representations can be constructed
explicitly, see \cite{CP2,CP3}. We will need from the construction
only the formulas for the spectra of the operators $\Phi^+(u)$.

The representation $W_r(a)$ has a basis $\{ v^{(r)}_i
\}_{i=0,\ldots,r}$, and $\Phi^\pm(u)$ acts on them as follows:
\begin{align}    \label{spec1}
\Phi^\pm(u) \cdot v^{(r)}_i &= q^{r-2i}
\frac{(1-uaq^{-r})(1-uaq^{r+2})}{(1-uaq^{r-2i+2})(1-uaq^{r-2i})}
v^{(r)}_i \\    \label{spec2}
&= q^{r-2i} \prod_{k=1}^r \frac{1-uaq^{r-2k}}{1-uaq^{r-2k+2}}
\prod_{j=1}^i \frac{1-uaq^{r-2j+4}}{1-uaq^{r-2j}} v^{(r)}_i
\end{align}
(we consider the right hand side as a power series in $u^{\pm 1}$).

Using these formulas we obtain the following expression for the
$q$--character of $W_r(a)$:
\begin{equation}    \label{wra}
\chi_q(W_r(a)) = \prod_{k=1}^r Y_{aq^{r-2k+1}} \left( \sum_{i=0}^r
\prod_{j=1}^i A_{aq^{r-2j+2}}^{-1} \right),
\end{equation}
where
\begin{equation}    \label{APhi}
A_a = Y_{aq} Y_{aq^{-1}} = q^2 \Phi^-(z^{-1} a^{-1}).
\end{equation}

Following \cite{CP2,CP3}, call the set $\Sigma_{a,r} = \{ aq^{r-2k+1}
\}_{k=1,\ldots,r}$ a $q$--segment of length $r$ and with center
$a$. Two $q$--segments are said to be in {\em special position} if
their union is a $q$--segment that properly contains each of them.

\begin{thm}[\cite{CP2,CP3}]    \label{classl2}
{\em The tensor product $W_{r_1}(b_1) \otimes \ldots \otimes
W_{r_m}(b_m)$ is irreducible if and only if none of the segments
$\Sigma_{r_i}(b_i)$ are in pairwise special position. Further, each
irreducible finite-dimensional representation of $\uqsl$ is isomorphic
to a tensor product of this form.}
\end{thm}

\thmref{classl2} suggests the following construction of the
irreducible representation $V(\prod_{i=1}^n (1-ua_i))$ (see
Chari-Pressley \cite{CP2,CP3}). One can break the set $\{ a_i \}$ in a
unique way into a union of segments $\Sigma_{r_i}(b_i)$, which are not
in pairwise special position. Then $V(P) \simeq W_{r_1}(b_1) \otimes
\ldots \otimes W_{r_m}(b_m)$.

\thmref{classl2} together with the multiplicative property of
$q$--characters imply that the $q$--character of any irreducible
finite-dimensional representation of $\uqsl$ is the product of the
$q$--characters given by formula \eqref{wra}. We conclude that the
$q$--character of the irreducible representation corresponding to
$P(u) = \prod_{i=1}^n (1-ua_i)$ equals
\begin{equation}    \label{charsl2}
Y_{a_1} \ldots Y_{a_n} \left( 1 + \sum_p M'_p \right),
\end{equation}
where each $M'_p$ is a monomial of the form $A_{c_1}^{-1} \ldots
A_{c_l}^{-1}$ with $c_j \in \bigcup a_i q^{2\Z}$.

Therefore we can associate to an irreducible finite-dimensional
representation $V$ of $\uqsl$ an oriented graph $\Gamma_V$, whose
vertices are labeled the monomials appearing in the $q$--character of
$V$. We connect the vertices corresponding to monomials $M_1$ and
$M_2$ by an arrow pointing towards $M_2$, if $M_2 = M_1 A_c^{-1}$ and
assign to this arrow the number $c$. It is clear that the graph
$\Gamma_V$ is connected.

\begin{definition}
Let $P(u)$ be a polynomial with the set of inverse roots $\{ a_i \}$
split in a unique way into a union of segments $\Sigma_{r_i}(b_i)$
that are not in pairwise special position. Then $P(u)$ is called {\em
irregular} if $\Sigma_{r_i}(b_i) \subset \Sigma_{r_j}(b_j)$ and
$\Sigma_{r_i}(b_iq^2) \subset \Sigma_{r_j}(b_j)$ at least for one pair
$i \neq j$. Otherwise, $P(u)$ is called {\em regular}.
\end{definition}

\begin{definition}
Given a ring of polynomials $\Z[x_\al^{\pm 1}]_{\al \in A}$, let us
call a monomial in this ring {\em dominant} if it has the form
$\prod_{k=1}^n x_{\al_k}$, i.e. it does not contain $x_\al^{-1}$.
\end{definition}

\begin{lem}    \label{domi}
The irreducible representation $V(P(u))$ contains dominant terms other
than the one corresponding to the highest weight vector if and only if
$P(u)$ is irredular.
\end{lem}

\begin{proof}
Let $\bigcup \Sigma_{r_i}(b_i)$ be the splitting of the set $\{ a_i
\}$ of inverse roots of $P(u)$ into a union of segments that are not
in pairwise special position. By \thmref{classl2}, $V(P(u)) \simeq
W_{r_1}(b_1) \otimes \ldots \otimes W_{r_m}(b_m)$. According to
formula \eqref{wra}, $\chi_q(W_r(a))$ is the sum of the monomials
\begin{equation}    \label{wra1}
\prod_{i=0}^{m-1} Y_{aq^{2i-r+1}} \prod_{j=m+1}^{r+1}
Y_{aq^{2j-r+1}}^{-1}, \quad \quad m=0,\ldots,r.
\end{equation}
Suppose that the product $\chi_q(W_{r_1}(b_1)) \ldots
\chi_q(W_{r_m}(b_m))$ contains a dominant term other the monomial
corresponding to the highest weight vector. Then this term is the
product of monomials $M_k$ from each $\chi_q(W_{r_k}(b_k))$, at least
one of which is not dominant. Without loss of generality we can assume
that $M_1$ is not dominant. But it is clear that if $M_1 M_2 \ldots
M_m$ is dominant, then so is $M_1 \wt{M}_2 \ldots \wt{M}_m$, where
$\wt{M}_k$ is the dominant term of $\chi_q(W_{r_k}(b_k))$, $\wt{M}_k =
\prod_{i=0}^{r_k-1} Y_{aq^{2i}}$. Formula \eqref{wra1} shows that for
$M_1 \wt{M}_2 \ldots \wt{M}_m$ to be dominant, the intersection
between $\Sigma_{r_1}(b_1q^2)$ and the union of the segments
$\Sigma_{r_j}(b_j), j\neq 1$, has to be non-empty. Since the segments
are not in pairwise special position by our assumption, this
immediately implies that $P(u)$ is irregular.

The converse statement is now also clear.
\end{proof}

For instance, if each $a_i$ has multiplicity $1$ (i.e., none of the
segments $\Sigma_{r_i}(b_i)$ is contained in another), then according
to \lemref{domi}, $\chi_q(V(\prod_{i=1}^n (1-ua_i)))$ has no dominant
terms other than $Y_{a_1} \ldots Y_{a_n}$.

\lemref{domi} implies the following result. Denote by
$\Z_+[x_\al]_{\al \in A}$ the subset of $\Z[x_\al^{\pm 1}]_{\al \in
A}$ consisting of all linear combinations of monomials in $x_\al^{\pm
1}$ with positive integral coefficients. Let $\Rep_{\on{reg}} \uqsl$
be the abelian subcategory of $\Rep \uqsl$ with the irreducible
objects $V(P(u))$ with regular polynomials $P(u)$. Denote by
$\on{Rep}_{\on{reg}} \uqsl$ the Grothendieck ring of this subcategory.

\begin{cor}    \label{regul}
Let ${\mc V} \in \on{Rep}_{\on{reg}} \uqsl$ be such that $\chi_q({\mc
V}) \in \Z_+[Y_a^{\pm 1}]_{a \in \C^\times}$. Then ${\mc V}$ is a
linear combination of irreducible representations of $\uqsl$ with
positive integral coefficients only.
\end{cor}

\begin{proof}
Suppose that there exist irreducible representations $V(P_i),
i=1,\ldots,n$, and $V(Q_j), j=1,\ldots,m$, where $P_i$'s and $Q_j$'s
are regular polynomials, such that
\begin{equation}    \label{empty}
\sum_{i=1}^n \chi_q(V(P_i)) = \sum_{j=1}^m \chi_q(V(Q_j)).
\end{equation}
The left hand side contains dominant terms $Y_{a^{(i)}_1} \ldots
Y_{a^{(i)}_{n_i}}$, where $P_i(u) = \prod_{k=1}^{n_i}
(1-ua^{(i)}_k)$. According to \lemref{domi}, the right hand side
contains such monomials if and only if each $P_i(u)$ equals some
$Q_j(u)$. Repeating this argument for the right hand side, we see that
each representation $V(P_i)$ is isomorphic to a unique representation
$V(Q_j)$, and therefore the relation \eqref{empty} is empty.
\end{proof}

Unfortunately, the statement of \corref{regul} is not true if ${\mc
V}$ is not assumed to lie in $\on{Rep}_{\on{reg}} \uqsl$, as the
following counter-example, due to E.~Mukhin, shows:
$$
\chi_q(W_1(a) \otimes W_2(aq)) + \chi_q(W_1(a) \otimes
W_2(aq^{-1})) - \chi_q(W_1(a)) =
$$
$$
Y_a^2 Y_{aq^2} + Y_a^2 Y_{aq^4}^{-1} + 2 Y_a Y_{aq^2}^{-1}
Y_{aq^4}^{-1} + Y_{aq^2}^{-2} Y_{aq^4}^{-1} + Y_a^2 Y_{aq^{-2}} + 2
Y_a Y_{aq^{-2}} Y_{aq^2}^{-1} + Y_{aq^{-2}} Y_{aq^2}^{-2} + Y_a^{-1}
Y_{aq^2}^{-2}.
$$
This is an element of $\Z_+[Y_a^{\pm 1}]_{a \in \C^\times}$, which
lies in the image of $\chi_q$, but is not the $q$--character of an
actual representation of $\uqsl$.

Finally, note that it is possible to write down a closed formula for
the $q$--characters of all irreducible $\uqsl$--modules. A similar
formula in the Yangian case (following Knight's definition of
character \cite{Kn}) was obtained by Chari and Pressley in \cite{CP5}.

\subsection{General case}    \label{gen case}

Now we can generalize some of the results of the previous section.

Let us set
\begin{equation}    \label{Aia}
A_{i,a} = q_i^2 k_i^{-1} \exp \left( - (q-q^{-1}) \sum_{n>0} h_{i,-n}
z^n a^n \right) = q_i^2 \Phi^-_i(z^{-1} a^{-1}), \quad \quad a \in
\C^\times.
\end{equation}
Clearly, $A_{i,a} \in \yy$. Using formula \eqref{connection}, we can
express $A_{i,a}$ in terms of $Y_{j,b}$'s:
\begin{equation}    \label{express}
A_{i,a} = Y_{i,aq_i} Y_{i,aq_i^{-1}} \prod_{j:I_{ji}=1} Y_{j,a}^{-1}
\prod_{j:I_{ji}=2} Y_{j,aq}^{-1} Y_{j,aq^{-1}}^{-1} \prod_{j:I_{ji}=3}
Y_{j,aq^2}^{-1} Y_{j,a}^{-1} Y_{j,aq^{-2}}^{-1}.
\end{equation}

\begin{prop}    \label{monom}
The $q$--character of the irreducible finite-dimensional
representation $V(\bp)$, where
\begin{equation}    \label{formbp}
P_i(u) = \prod_{k=1}^{n_i} (1-ua^{(i)}_k), \quad \quad i \in I,
\end{equation}
equals
\begin{equation}    \label{mp}
\prod_{i\in I} \prod_{k=1}^{n_i} Y_{i,a^{(i)}_k} \left( 1 + \sum M'_p
\right),
\end{equation}
where each $M'_p$ is a monomial in $A_{j,c}^{\pm 1}$.
\end{prop}

In order to prove this, we consider the $R$--matrices
$$
R_{V,W}(z) = (\pi_{V(z)} \otimes \pi_W)(\R) = \pi_W(\R) \in \on{End}(V
\otimes W)[[z]].
$$
Let us recall the general result about the $R$--matrices (see
\cite{IFR} and \cite{EFK}, Prop. 9.5.3), which follows from their
crossing symmetry property.

\begin{prop}    \label{factor}
For any pair of irreducible representations $V=V(\bp), W=V(\bq)$ of
$\uqg$,
$$
R_{V,W}(z) = f_{V,W}(z) \ovl{R}_{V,W}(z),
$$
where the matrix elements of $\ovl{R}_{V,W}(z)$ are the
$z$--expansions of rational functions in $z$ that are regular at
$z=0$, with $\ovl{R}_{V,W}(z) \cdot v_{\bp} \otimes v_{\bq} = v_{\bp}
\otimes v_{\bq}$, and $f_{V,W}(z)$ can be represented in the form
$$
f_{V,W}(z) = q^{-(\la,\mu)} \prod_{n=1}^\infty
\rho_{V,W}(zq^{2r^\vee h^\vee n}),
$$
where $\la$ and $\mu$ are the degrees of $\bp$ and $\bq$,
respectively, and $\rho(z)$ is the expansion of a rational function.
\end{prop}

Now observe that 
\begin{equation}     \label{eigenvalue}
t_V(z) \cdot v_{\bq} = {\mathbf h}_q(t_V(z)) \cdot v_{\bq} = \chi_q(V)
\cdot v_{\bq},
\end{equation}
so $v_{\bq}$ is an eigenvector of $t_V(z)$. According to
\thmref{mainthm},(1), $\chi_q(V)$ is a polynomial in $Y_{i,a}^{\pm
1}$. Recall that $Y_{i,a}$ is a power series in $z$, whose
coefficients are polynomials in the generators $k_j,h_{j,n}, j \in I,
n<0$, and rational functions in $q^n$. We have:
$$
Y_{i,a}^{\pm 1} \cdot v_{\bq} = Y_i^{\bq}(za)^{\pm
1} v_{\bq},
$$
where $Y_i^{\bq}(z) \in \C[[z]]$. Thus, the eigenvalue of $t_V(z)$ on
$v_{\bq}$ equals $\chi_q(V)$, in which we substitute each $Y_{i,a}$ by
$Y_i^{\bq}(za)$. Each $Y_i^{\bq}(za)$ can be found if we solve the
equations \eqref{express} for $i \in I$, since we know that the value
of $A_{i,a} = q_i^2 \Phi^-_i(z^{-1}a^{-1})$ on $v_{\bq}$ equals
$q_i^{2+\deg Q_i}
\dfrac{Q_i(z^{-1}a^{-1}q_i^{-1})}{Q_i(z^{-1}a^{-1}q_i)}$.

On the other hand, let us choose bases of generalized eigenvectors of
$\Phi^\pm(u)$ in $V$ and $V(\bq)$, so that the latter includes vector
$v_{\bq}$. The eigenvalue of $t_V(z)$ on the highest weight vector
$v_{\bq} \in V(\bq) = W$ equals the sum of the diagonal entries of the
$R$--matrix $R_{V,W}(z)$ written in this basis, which correspond to
the vectors $v \otimes v_{\bq}$, where $v$ is a basis vector of $V$
(up to $q^\rho$). The proof of \thmref{mainthm} shows that these
diagonal entries are in one-to-one correspondence with the monomials
appearing in the $q$--character $\chi_q(V)$, namely, the diagonal
entry corresponding to the monomial $Y_{i_1,a_1} \ldots Y_{i_k,a_k}$
equals $Y_{i_1}^{\bq}(za_1) \ldots Y_{i_k}^{\bq}(za_k)$. This
observation allows us to compute $f_{V,W}(z)$ explicitly.

\subsection{Computation of $f_{V,W}(z)$}
Consider the diagonal matrix element $f_{V,W}(z)$ of $R_{V,W}(z)$
corresponding to the vector $v_{\bp} \otimes v_{\bq}$. It has been
computed in some cases (see, e.g., \cite{FR:crit,AK,EFK}). The
following proposition gives a general formula for an arbitrary pair of
irreducible representations of $\uqg$.

\begin{prop}
Let $\bp = (P_i)_{i\in I}$, where
$$
P_i(u) = \prod_{k=1}^{n_i} (1-ua^{(i)}_k), \quad \quad i \in I.
$$
Then
$$
f_{V(\bp),V(\bq)}(z) = q^{-(\la,\mu)} \prod_{i\in I} \prod_{k=1}^{n_i}
Y_i^{\bq}(za^{(i)}_k).
$$
\end{prop}

To find $Y_i^{\bq}(z)$ explicitly, note first that if $Q_j(u) =
\prod_{l=1}^{m_j} (1-ub^{(j)}_l), j \in I$, then
\begin{equation}    \label{product}
Y_i^{\bq}(z) = \prod_{j\in I} \prod_{l=1}^{m_j}
Y_i^{(j)}(z/b^{(j)}_l),
\end{equation}
where $Y_i^{(j)}(z/b)$ corresponds to $\bq = \bp^{(i)}_b$, where
$(\bp^{(j)}_b)_i = (1-ub)$, if $i=j$ and $1$, if $i\neq j$ (this is
the case when $W = V_{\omega_j}(b)$).

In the same way as in the proof of \thmref{mainthm} we find that the
eigenvalue of $h_{i,-n}, n>0$, on the highest weight vector of
$V_{\omega_j}(b)$ equals $\delta_{i,j} (q_j^n-q_j^{-n})
b^{-n}/n(q-q^{-1})$. Using formula \eqref{connection}, we find that
the eigenvalue of $\wt{h}_{i,-n}$ on this vector equals
\beq    \label{inverse}
\frac{b^{-n}(q_j^n-q_j^{-n})}{n(q-q^{-1})} \wt{C}_{ji}(q^n).
\end{equation}
The matrix $\wt{C}(x) = C(x)^{-1}$ is known explicitly for $\g$ of
classical types (see Appendix C of \cite{FR:simple}), and for an
arbitrary $\g$ the determinant of $C(x)$ is known (see
\cite{BP}). These results imply that $\wt{C}_{ij}(x)$ equals
$\ovl{C}_{ij}(x)/(1-x^{2r^\vee h^\vee})$, where $\ovl{C}_{ij}(x)$ is a
polynomial in $x$ with integral coefficients. When we substitute it
into formula \eqref{inverse} and then into formula \eqref{Yia}, we
obtain that each monomial $\pm x^p$ appearing in $\ovl{C}_{ij}(x)$
contributes the factor
\begin{equation}    \label{fact}
(zq^pq_i;q^{2r^\vee h^\vee})_\infty^{\pm 1} (zq^pq_i^{-1};q^{2r^\vee
h^\vee})_\infty^{\mp 1}
\end{equation}
to $Y_i^{(j)}(z)$, where
$$
(a;b)_\infty = \prod_{n=1}^\infty (1-ab^n).
$$

Thus, for any irreducible representations $V, W$, the function
$f_{V,W}(z)$ equals $q^{-(\la,\mu)}$ times the product of the factors
$(zq^n;q^{2r^\vee h^\vee})_\infty^{\pm 1}$, where $n \in \Z$.  This
statement is stronger than the corresponding statement of
\propref{factor} in that we claim that the zeroes and poles of
$\rho(x)$ are necessarily integral powers of $q$. This result can be
used to gain insights into the structure of the poles of the
$R$--matrices.

Following Akasaka and Kashiwara \cite{AK}, let us denote by
$d_{V,W}(z)$ the denominator of $\ovl{R}_{V,W}(z)$, i.e., the
polynomial in $z$ of smallest degree, such that $d_{V,W}(z)
\ovl{R}_{V,W}(z)$ has no poles. We normalize it so that its constant
term is equal to $1$. Using the crossing symmetry of the $R$--matrices
one can show (see \cite{AK}, Prop.~A.1) that $f_{V,W}(z)$ satisfies
\begin{equation}    \label{diffeq}
f_{V,W}(z)f_{^* V,W}(z) = c \frac{d_{V,W}(z)}{d_{W,^* V}(z^{-1})},
\end{equation}
where $c = c' z^n, c' \in \C^\times, n \in \Z$.

As we have shown above, the left hand side is a rational function,
whose zeroes and poles are powers of $q$. Therefore if $y$ is a pole
of $R_{V,W}(z)$ (i.e., a root of $d_{V,W}(z)$) that is not a power of
$q$, then $y^{-1}$ has to be a pole of $R_{W,^* V}(z)$. It is unclear
to us at the moment whether one can use this kind of argument to prove
that the poles of $R_{V,W}(z)$ are integral powers of $q$ as
conjectured in \cite{AK}.

\begin{rem}    \label{comp}
Here we have used the $q$--characters to obtain information about the
$R$--matrices. Conversely, we can use the information about the
$R$--matrices to gain insights into the structure of the
$q$--characters. In fact, $\chi_q(W)$ can be read off the diagonal
entries of the $R$--matrices $R_{W,V_{\omega_i}(1)}(z), i \in I$,
corresponding to the highest weight vectors in $V_{\omega_i}(1)$.\qed
\end{rem}

\subsection{Proof of \propref{monom}.} Consider the diagonal entry of
the $R$--matrix corresponding to a monomial
$$
\prod_{i\in I} \prod_{k=1}^{n_i} Y_{i,a^{(i)}_k} \cdot \prod_{k=1}^m
Y_{i_k,a_k}^{\ep_k}
$$
occurring in $\chi_q(V)$. Then this entry of the $R$--matrix equals
$$
f_{V(\bp),V(\bq)}(z) \cdot \prod_{k=1}^m
Y_{i_k}^{\bq}(za_k)^{\ep_k}.
$$
\propref{factor} then implies that $\prod_{k=1}^m
Y_{i_k}^{\bq}(za_k)^{\ep_k}$ is a rational function in $z$ for all
$\bq$. According to formula \eqref{product}, this is equivalent to its
being a rational function when $\bq = \bp^{(j)}_b, \forall j \in I, b
\in \C^\times$. The following lemma shows that in that case
$\prod_{k=1}^m Y_{i_k,a_k}^{\ep_k}$ can be written as the product of
$A_{i,c}^{\pm 1}$. This proves \propref{monom}.

\begin{lem} Suppose that a monomial $M = \prod_{k=1}^m
Y_{i_k,a_k}^{\ep_k}$ is such that
$\prod_{k=1}^m Y_{i_k}^{(j)}(za_k)^{\ep_k}$ is a rational function for
all $j \in I$. Then $M$ can be written as a product of $A_{i,c}^{\pm
1}$.
\end{lem}

\begin{proof}
We know that each $Y_{i_k}^{(j)}(za_k)^{\ep_k}$ is a product of the
factors \eqref{fact}. Therefore, without loss of generality, we can
assume that each $a_k = a q^{l_k}$ for some $a \in \C^\times $ and
$l_k \in \Z$. According to formula \eqref{Yia}, the above condition on
$M$ then means that the eigenvalue of
$$
a^n \sum_{k=1}^m \ep_k \wt{h}_{i_k,-n} q^{l_k n}
$$
on $v_{\bp^{(j)}_1}$ equals $a^n \beta(q^n)/n$, where $\beta(x)$ is a
polynomial in $x^{\pm 1}$ with integral coefficients. In other words,
there exist polynomials $\beta(x)$ and $\ga_i(x), i \in I$, such that
the eigenvalue of
$$
(q-q^{-1}) \sum_{i \in I} \wt{h}_{i,-n} \ga_i(q^n)
$$
is equal to $\beta(q^n)/n$. But
$$
\sum_{i \in I} \wt{h}_{i,-n} \ga_i(q^n) = (q-q^{-1}) \sum_{i,j \in I}
h_{j,-n} \wt{C}_{ji}(q^n) \ga_i(q^n),
$$
and since the eigenvalue of $h_{i,-n}$ on $v_{\bp^{(j)}_1}$ is
$\delta_{i,j}[n]_{q_i}/n$, we obtain that
\begin{equation}    \label{poli}
\sum_{i \in I} \wt{C}_{ji}(x) \ga_i(x) (x^{r_j} - x^{-r_j})
\end{equation}
is a polynomial in $x^{\pm 1}$. The lemma says that for \eqref{poli}
to be a polynomial for all $j \in I$, $\ga_i(x)$ must have the form
$$
\ga_i(x) = \sum_{k\in I} C_{ik}(x) R_k(x), \quad \quad i \in I,
$$
where $R_k(x)$ are some polynomials in $x^{\pm 1}$.

It is clear that if we allow $\ga_i(x)$ to be a rational function, and
we want \eqref{poli} to be equal to a polynomial $R'_j(x)$ for all $j
\in I$, then $\ga_i(x)$ can be represented in the form
\begin{equation}    \label{poli1}
\ga_i(x) = \sum_{k\in I} C_{ik}(x) \frac{R'_k(x)}{x^{r_k}-x^{-r_k}},
\quad \quad i \in I.
\end{equation}
It remains to show that \eqref{poli1} is a polynomial for all $i \in
I$ if and only if each $R'_k(x)$ is divisible by $x^{r_k}-x^{-r_k}$.

Since $\det C(\pm 1) \neq 0$, $R'_i(x)$ is divisible by $x-x^{-1}$. This
proves the result for simply laced $\g$, for which $r_k=1, \forall k
\in I$. For non-simply laced $\g$, we can now replace
$1/(x^{r_k}-x^{-r_k})$ in formula \eqref{poli1} with
$(x-x^{-1})/(x^{r_k}-x^{-r_k})$. After that the result is easy to
establish by inspection.
\end{proof}

\subsection{The conjecture}
The following conjecture is motivated by explicit examples of the
$q$--characters (see Sects.~\ref{casesl2},\ref{ex}).

\begin{conj}    \label{genconj}
{\em The $q$--character of the irreducible finite-dimensional
representation $V(\bp)$, where $\bp$ is given by formula
\eqref{formbp}, can be represented in the form \eqref{mp} where each
$M'_p$ is a monomial in $A_{j,c}^{-1}, c \in \C^\times$.}
\end{conj}

\conjref{genconj} holds for $U_q \sun$. In this case the
$q$--characters of $V_{\omega_i}(a)$ can be found (see \remref{comp})
from the explicit formulas for the $R$--matrices
$R_{V_{\omega_i}(z),V_{\omega_j}(w)}, i,j=1,\ldots,N-1$, which are
known in the literature (see \cite{DatO,AK}). For other classical $\g$
and $\g=G_2$, \conjref{genconj} can probably also be derived from a
case by case analysis based on what is known about the structure of
the $R$--matrices of the fundamental representations and the
decompositions of their tensor products. The expected formulas for the
$q$--characters of the fundamental representations for these algebras
are given in \cite{FR:crit,FR:simple} (see also \secref{ex} below) for
classical $\g$, and \cite{Ko,BP} for $\g=G_2$.

Now we give two corollaries to \conjref{genconj}.

\begin{cor}    \label{cor1}
\hfill

(1) $\chi_q(V_{\omega_i}(a))$ equals $Y_{i,a}(1 + \sum_p M'_p)$, where
each $M'_p$ is a monomial in $A_{j,aq^n}^{-1}$, $n \in \Z$, and
$Y_{i,a}$ is the only dominant monomial in $\chi_q(V_{\omega_i}(a))$.

(2) The monomials $M'_p$ occurring in $\chi_q(V(\bp))$, where $\bp$ is
given by formula \eqref{formbp}, are products of $A_{j,c}^{-1}$, where
$c \in \bigcup a^{(i)}_k q^{\Z}$.
\end{cor}

\begin{proof}
Let us consider $\uqg$ as a module over $\C(q)$. It then has an
algebra automorphism that sends $q$ to $q^{-1}$, $k_i$ to $k_i^{-1}$
and leaving the generators $x^\pm_i$ unchanged. If we apply this
automorphism to $\R$, we obtain $(S \otimes \on{id})(\R)$, where $S$
is the antipode. This means that for any representation $V$ of $\uqg$,
if we replace in $\chi_q(V)$ each $Y_{i,a}$ by $Y_{i,\ovl{a}}^{-1}$,
where $\ovl{a}$ is obtained from $a$ by replacing $q$ by $q^{-1}$,
then we obtain $\chi_q(V^*)$.

On the other hand, \cite{CP4}, Proposition 5.1(b) implies that
$V_{\omega_i}(a)^* \simeq V_{\omega_{\ovl{i}}}(ap^{-1})$, where $p =
q^{r^\vee h^\vee}$ and $\omega_{\ovl{i}} = -w_0(\omega_i)$, $w_0$
being the longest element of the Weyl group of $\g$.

Now, according to \propref{monom}, we can write
$$
\chi_q(V_{\omega_i}(a)) = Y_{i,a}(1 + \sum_r M'_r), \quad \quad
\chi_q(V_{\omega_{\ovl{i}}}(a)) = Y_{\ovl{i},a}(1 + \sum_r M''_r),
$$
where each $M'_r, M''_r$ is a monomial in $A_{j,c}^{\pm 1}$, and
\begin{equation}    \label{YY}
Y_{i,a}^{-1} Y_{\ovl{i},ap^{-1}}^{-1} = \ovl{M}'_r M''_r.
\end{equation}
Here $\ovl{M}'_r$ is obtained from $M'_r$ be replacing each $A_{j,c}$
by $A_{j,\ovl{c}}$.

But \conjref{genconj} tells us that both $M'_r$ and $M''_r$ are
monomials in $A_{j,c}^{-1}$ only. Therefore they can only be monomials
in $A_{j,q^n}^{-1}, n \in \Z$, for otherwise formula \eqref{YY} can
not hold. This implies that $\chi_q(V_{\omega_i}(a))$ equals
$Y_{i,a}(1 + \sum_p M'_p)$, where each $M'_p$ is a monomial in
$A_{j,aq^n}^{-1}, n \in \Z$.

Let us now prove that the only dominant term in
$\chi_q(V_{\omega_i}(a_i))$ is $Y_{i,a}$.

In the same way as in the proof of \propref{monom}, we can show that
this statement is equivalent to the following. If $R_k(x), k \in I$,
are polynomials in $x^{\pm 1}$ with non-negative integral
coefficients, such that
$$
- \sum_{k\in I} C_{jk}(x) R_k(x) + \delta_{ij}
$$
is a polynomial with non-negative integral coefficients for all $j\in
I$, then $R_k(x)=0, \forall k \in I$ (the matrix $C_{ij}(x)$ is given
in \secref{cartan}). Suppose that this is not so. Let $l_k$ and $h_k$
be the lowest and highest degrees of those $R_k(x)$ that are
non-zero. Choose the smallest, $l_s$, among $l_k$'s, and the largest,
$h_t$, among $h_k$'s. If there are several $R_k(x)$ with the same
lowest (resp., highest) degree, we pick the $s$ (resp. $t$) that
corresponds to the largest value of $r_k$ (i.e., to the longer
root). Then we obtain that $- \sum_{k\in I} C_{sk}(x) R_k(x)$ contains
the monomial $x^{l_s-r_s}$ with a negative coefficient. This monomial
can only be compensated by $\delta_{is}$. But then $l_s=r_s$. Applying
the same argument to $h_t$, we obtain that $h_t=-r_t$, and so
$h_t<l_s$, which is a contradiction. This completes the proof of part
(1) of the corollary.

Part (2) follows from the fact that all other representations can be
obtained as subfactors of the tensor products of the fundamental
representations.
\end{proof}

\begin{cor}    \label{not}
The tensor product $V_{\omega_1}(a_1) \otimes \ldots \otimes
V_{\omega_n}(a_n)$ is irreducible if $a_j/a_k \not{\hspace*{-1mm}\in}
q^{\Z}$, $\forall i\neq j$.
\end{cor}

\begin{proof}
If $V_{\omega_1}(a_1) \otimes \ldots \otimes V_{\omega_n}(a_n)$ is
reducible, then $\chi_q(V_{\omega_1}(a_1)) \ldots
\chi_q(V_{\omega_n}(a_n))$ should contain a dominant term other than
the product of the highest weight terms. But for that to happen, for
some $j$ and $k$, there have to be cancellations between some
$Y_{p,a_jq^{n}}^{-1}$ appearing in $\chi_q(V_{\omega_j}(a_j))$ and
some $Y_{r,a_kq^{m}}$ appearing in $\chi_q(V_{\omega_k}(a_k))$. These
cancellations may only occur if $a_j/a_k \in q^{\Z}$.
\end{proof}

\begin{rem} \corref{not} has been conjectured earlier
by Akasaka and Kashiwara \cite{AK}. Furthermore, they conjectured that
$V_{\omega_i}(a) \otimes V_{\omega_j}(b)$ is reducible only when
$b=aq^n$, where $n \in \Z, |n| \leq r^\vee h^\vee$. This can also be
derived from \conjref{genconj}, because it is easy to see that
$\chi_q(V_{\omega_i}(a))$ is a linear combination of monomials in
$Y_{k,q^n}^{\pm 1}$, where $n \in \Z, 0 \leq n \leq r^\vee h^\vee$. We
thank Kashiwara for a discussion of these conjectures.\qed
\end{rem}

\section{Combinatorics of $q$--characters}

\subsection{Interpretation in terms of the joint spectra of
$\Phi_i^\pm(u)$}

For $i \in I, a \in \C^\times$, we define special polynomials
$R^{i,a}_j(u)$ and $Q^{i,a}_j(u)$. Consider the matrix element
$C_{ji}(q)$ of the matrix $C(q)$ defined in \secref{cartan}. This is a
combination of powers of $q$ with coefficients $\pm 1$. We define
$R^{i,a}_j(u)$ (resp., $Q^{i,a}_j(u)$) as the polynomial in $u$ with
constant term $1$, whose zeroes are $a$ times the powers of $q$
appearing in $C_{ji}(q)$ with coefficient $1$ (resp.,
$-1$). Explicitly, we have:
\begin{align*}
R^{i,a}_i(u) &= (1-uaq^{r_i})(1-uaq^{-r_i}), \quad \quad Q^{i,a}_i(u)
= 1,\\
R^{i,a}_j(u) &= 1, \quad \quad i\neq j,\\
Q^{i,a}_j &= \left\{ \begin{array}{ll}
1, \quad  & I_{ji}=0,\\
1-ua, \quad & I_{ji}=1,\\
(1-uaq)(1-uaq^{-1}), \quad & I_{ji}=2,\\
(1-uaq^2)(1-ua)(1-uaq^{-2}), \quad & I_{ji}=3
\end{array} \right. \quad \quad i\neq j.
\end{align*}
Finally, let us set
$$
P^{i,a}_j(u) = Q^{i,a}_j(u) R^{i,a}_j(u)^{-1}.
$$

Now it is clear how to interpret \propref{monom} and \conjref{genconj}
in terms of the eigenvalues of $\Phi^\pm_i(u)$. For instance,
\propref{monom} means that the eigenvalues of $\Phi^\pm_i(u)$ on the
module $V(\bp)$ have the form
$$
\frac{P_i(uq_i^{-1})}{P_i(uq_i)} \frac{\prod_{k=1}^m
P^{i,c_k}_{j_k}(uq_i^{-1})^{\pm 1}}{\prod_{k=1}^m
P^{i,c_k}_{j_k}(uq_i)^{\pm 1}},
$$
for some $c_1,\ldots,c_m \in \C^\times$ (up to obvious overall power
of $q$ factors). Explicit calculation shows that
\begin{equation}    \label{piaj}
\frac{P^{i,a}_j(uq_i^{-1})}{P^{i,a}_j(uq_i)} =
\frac{1-uaq^{(\al_i,\al_j)}}{1-uaq^{-(\al_i,\al_j)}}.
\end{equation}

There is an interesting connection between the eigenvalues of
$\Phi^\pm_i(u)$ (and hence the $q$--characters) and the action of
the generators $x_{i,n}^-$ on the finite-dimensional $\uqg$--modules.

Let us consider again the $\uqsl$--module $W_r(a)$. The action of the
generators $x^-_n$ in the basis $\{ v^{(r)}_i \}$ is given by a very
simple formula:
$$
x_n^- \cdot v^{(r)}_i = (aq^{r-2i})^n \; v^{(r)}_{i+1} = (aq^{r-2i})^n
\; x^-_0 \cdot v^{(r)}_i.
$$

Now suppose, more generally, that $v$ is a vector in a
finite-dimensional $\uqg$--module $V$, which is an eigenvector of
$\Phi^\pm_i(u), \forall i \in I$ with the eigenvalues $\Psi^\pm_i(u)$,
and that we have:
\begin{equation}    \label{xjn}
x^-_{j,n} \cdot v = \gamma^n \; x^-_{j,0} \cdot v, \quad \quad n \in \Z,
\end{equation}
for some $j \in I$ and $\gamma \in \C^\times$. In that case, using the
commutation relations between $h_{i,m}$ and $x^-_{j,n}$, we obtain
that $x^-_{j,0} \cdot v$ is also an eigenvector of $\Phi^\pm_i(u),
\forall i \in I$, with the eigenvalues
$$
\Psi^\pm_i(u) q^{(\al_i,\al_j)} \frac{1-u\gamma
q^{(\al_i,\al_j)}}{1-u\gamma q^{-(\al_i,\al_j)}}.
$$
According to formula \eqref{piaj}, we can rewrite this as
$$
\Psi^\pm_i(u) q^{(\al_i,\al_j)}
\frac{P^{i,\gamma}_j(uq_i^{-1})}{P^{i,\gamma}_j(uq_i)}.
$$
Therefore if $M$ is the monomial in $Y_{i,a}^{\pm 1}$ that corresponds
to vector $v$ in $\chi_q(V)$, then the monomial corresponding to the
vector $x^-_{j,0} \cdot v$ is $M \cdot A_{j,\gamma}^{-1}$.

This observation clarifies the statements of \propref{monom} and
\conjref{genconj} in the case when $V$ is spanned by vectors obtained
by the action of $x^-_{j,0}$ on the highest weight vector. Empirical
evidence suggests that in that case the action of the operators
$x^-_{n,j}$ indeed has the form \eqref{xjn}. Certainly, not all
representations have this property: for one thing, it means that the
restriction of $V$ to $\uqg$ is irreducible, which is not always the
case. Even in the case of $\uqsl$, the modules $W_r(a), r>0$, seem to
be the only representations with this property. Still, such
representations apparently exist in general, and it is plausible that
one can use them to prove \conjref{genconj}.

In the following two sections we discuss the combinatorial structure of
$q$--characters assuming that \conjref{genconj} is true.

\subsection{Reconstructing the $q$--character from the highest weight}
\label{gencase}

Recall the restriction property of the $q$--characters from
\thmref{mainthm},(3): if we apply the homomorphism $\on{res}_{\{i\}}$
to $\chi_q(V)$, i.e., replace all $Y_{j,b}, j \not{\hspace*{-1mm}\in}
I$ in $\chi_q(V)$ by $1$, we obtain the $q_i$--character of the
semi-simplification of the restriction of $V$ to $\uqg_{\{i\}} \simeq
U_{q_i} \widehat{\sw_2}$.

Given an irreducible representation $W$ of $U_{q_i}
\widehat{\sw_2}$ we write its $q$--character in the form
\eqref{charsl2}. We then replace each $A_c^{-1}$ by $A_{i,c}^{-1}$ and
denote the resulting element of $\Z[Y_{j,a}^{\pm 1}]_{j\in I}$ by
$\chi_q^{(i)}(W)$. It is clear that $\chi_q^{(i)}$ extends linearly to
a homomorphism $\on{Rep} U_{q_i} \widehat{\sw_2} \arr
\Z[Y_{i,a}^{\pm}]$, and its image equals $\Z[Y_{i,b} + Y_{i,b}
A_{i,bq_i}^{-1}]_{b \in \C^\times}$.

In view of \thmref{mainthm},(3) and \conjref{genconj} it is natural
to expect that for any $V \in \on{Rep} \uqg$,
$$
\chi_q(V) \in \bigcap_{i\in I} {\mc R}_i,
$$
where
\begin{equation}
{\mc R}_i = \Z[Y_{j,a_j}]_{j\neq i; a_j \in \C^\times} \otimes
\Z[Y_{i,b} + Y_{i,b} A_{i,bq_i}^{-1}]_{b \in \C^\times}.
\end{equation}
If $V$ is an actual representation, then $\chi_q(V) \in \bigcap_{i\in
I} {\mc R}_{i,+}$, where ${\mc R}_{i,+} = {\mc R}_i \cap
\Z_+[Y_{j,a}^{\pm}]_{j\in I}$.

It is natural to conjecture that {\em any} element of $\bigcap_{i\in
I} {\mc R}_{i,+}$ equals $\chi_q(V)$ for some representation $V$ of
$\uqg$. In other words, if an element of $\Z[Y_{i,a_i}]_{i \in I; a_i
\in \C^\times}$ has good restrictions, i.e., it restricts to
$q$--characters of $\uqg_{\{i\}}$ for each $i \in I$, then it is
necessarily a $q$--character of $\uqg$. This conjecture is essentially
equivalent to \conjref{main} that we state in \secref{screen}. Some
evidence for \conjref{main} coming from the theory of $\W$--algebras
is presented in \secref{defw}.

If this conjecture is true, then the $q$--character of the irreducible
representation $V(\bp)$, where $P_i(u) = \prod_{k=1}^{n_i}
(1-ua^{(i)}_k)$, can be constructed combinatorially as follows.

We start with the monomial $M = \prod_{i\in I} \prod_{k=1}^{n_i}
Y_{i,a^{(i)}_k}$ corresponding to the highest weight vector and try to
reconstruct the $q$--character of $V(\bp)$ from it. We know that any
element of ${\mc R}_{i,+}$ that contains a monomial of the form $N
\cdot \prod_{k=1}^{n_i} Y_{i,a^{(i)}_k}$, where $N \in \Z[Y_{j,a}^{\pm
1}]_{j\neq i}$, also contains $N \cdot
\chi_q^{(i)}(V(\prod_{k=1}^{n_i} (1-ua^{(i)}_k))$. Thus, the first
step is to replace the monomial $M$ with the product
$$
\prod_{i\in I} \chi_q^{(i)} \left( V(\prod_{k=1}^{n_i} (1-ua^{(i)}_k)
\right).
$$
Next, we continue by induction following the rule: whenever $\chi_q(V)$
contains a monomial $\prod_{k=1}^m Y_{i,b^{(i)}_k} N$, where $N \in
\Z[Y_{j,a}^{\pm 1}]_{j\neq i}$, it should be part of a combination
$$\chi_q^{(i)} \left( V(\prod_{k=1}^m (1-ub^{(i)}_k) \right) N.$$ Since
we know that the representation $V(\bp)$ is finite-dimensional, the
inductive process should stop after a finite number of steps. The
result should be $\chi_q(V(\bp))$. At each step we use only our
knowledge of the $q$--characters of $\uqsl$.

\begin{rem}
It is clear from the proof of \thmref{mainthm} that if ${\mc V} \in
\on{Rep} \uqg$ is a linear combination of irreducible representations
of $\uqg$ with positive integral coefficients, then $\chi_q({\mc V})
\in \Z_+[Y_{i,a_i}^{\pm 1}]_{i \in I; a \in \C^\times}$. The
counter-example given in \secref{casesl2} shows that the converse is
not true. However, it is possible that a weaker version still
holds.\qed
\end{rem}

\subsection{The graph $\Gamma_V$}

Now we attach to each irreducible finite-dimensional representation
$V$ of $\uqg$, an oriented colored graph $\Gamma_V$. Its vertices are
labeled by the monomials appearing in the $q$--character of
$V$.

Denote the monomial $\prod_{i\in I} \prod_{k=1}^{n_i} Y_{i,a^{(i)}_k}
M'_p$ by $M_p$. Two vertices corresponding to monomials $M_1$ and
$M_2$ are connected by an arrow pointing towards $M_2$, if $M_2 = M_1
A_{i,c}^{-1}$ for some $c$, {\em and} in
$\beta_{\{i\}}(\chi_q(V))$, the monomial $\beta_{\{i\}}(M_2)$
does not correspond to a highest weight vector of the
semi-simplification of $V|_{\uqg_{\{i\}}}$. We then assign to this
arrow the color $i$ and the number $c$.

It follows from the construction that if we erase in $\Gamma_V$ all
arrows but those of color $i$, we obtain the graph $\Gamma^{(i)}_V$ of
the semi-simplification of $V|_{\uqg_{\{i\}}}$. Furthermore, the
graphs $\Gamma_V$ are compatible with restrictions to all of its
quantum affine subalgebras: if $J$ is a subset of $I$, then the graph
of the semi-simplification of the restriction of $V$ to $\uqg_J$ can
be obtained from the graph $\Gamma_V$ by erasing the arrows of colors
$i \not{\hspace*{-1mm}\in} J$.

\conjref{genconj} implies that the graphs of the fundamental
representations $V_{\omega_i}(a)$ are connected, and moreover, there
exists an oriented path from the vertex corresponding to the highest
monomial $Y_{i,a}$ to any other vertex. Indeed, we have shown in the
proof of \corref{not} that $\chi_q(V_{\omega_i}(a))$ is the sum of
monomials of the form $Y_{i,a} \prod_k A_{i_k,q^{n_k}}^{-1}$, and the
only dominant term in $\chi_q(V_{\omega_i}(a))$ is $Y_{i,a}$. Now let
$M$ be an arbitrary monomial in $\chi_q(V_{\omega_i}(a))$ different
from $Y_{i,a}$. Then it is not dominant. Suppose that there are no
incoming arrows for this vertex. But then $\beta_{\{i\}}(M)$
corresponds to a highest weight vector in the semi-simplification of
$V|_{\uqg_{\{i\}}}$ for each $i \in I$. Hence $M$ is dominant, which
is a contradiction.

Now we see that $M$ has an incoming arrow, we can move up along this
arrow. The weight of the monomial on the other end of this arrow
equals that of $M$ plus $\al_i$. Continuing by induction, we arrive at
the highest weight monomial $Y_{i,a}$. Thus, we find an oriented path
from the vertex corresponding to the highest monomial $Y_{i,a}$ to the
vertex corresponding to a monomial $M$.

We conjecture that the graph $\Gamma_V$ is connected for any
irreducible representation $V$ of $\uqg$. It this is true, then the
reason for irreducible $\uqg$--modules being reducible when restricted
to $\uqgg$ essentially lies in that happening already for each of the
$U_{q_i} \widehat{\sw}_2$ subalgebras of $\uqg$.

\begin{rem}
The graph $\Gamma_V$ is similar to the crystal graph of
$V$. However, there is an important difference: while the arrows of a
crystal graph are labeled by the simple roots of $\G$, i.e., from $0$
to $\el$, the arrows of $\Ga_V$ are labeled by the simple roots of
$\g$, i.e., from $1$ to $\el$, and there is also a number attached to
each row. The crystal graph is designed so that it respects the
subalgebras $U_{q_i} \sw_2$, corresponding to the Drinfeld-Jumbo
realization of $\uqg$, while $\Ga_V$ respects the affine subalgebras
$U_{q_i} \widehat{\sw}_2$ in the Drinfeld ``new'' realization. It
would be interesting to develop a theory analogous to the theory of
crystal basis \cite{Ka} for the graphs $\Ga_V$.\qed
\end{rem}

\subsection{Examples}    \label{ex}

Here we give examples of $q$--characters and graphs associated to the
first fundamental representation of $\uqg$, where $\g$ is of classical
type. These $q$--characters can be obtained (see \remref{comp}) from
the explicit formulas for the $R$--matrices
$R_{V_{\omega_1}(z),V_{\omega_j}(w)}$ that are known in the
literature: see \cite{DatO,AK} for $A_\el$, \cite{AK} for $C_\el$, and
\cite{DavO} for $B_\el$ and $D_\el$. They also agree with the formulas
for the eigenvalues of the transfer-matrices \cite{Re1,Re2,BR,KS1}.

In all cases,
$$
\chi_q(V_{\omega_1}(a)) = \sum_{i \in J} \La_{i,a}.
$$

Additional examples can be found in \cite{FR:crit,FR:simple,Ko}.

\subsubsection{The $A_\el$ series}    \label{exa}

$J = \{ 1,\ldots,\el+1 \}$.

$$
\La_{i,a} = Y_{i,aq^{i-1}} Y_{i-1,aq^{i}}^{-1}, \quad
\quad i=1,\ldots,\el+1.
$$

Equivalently,
\begin{align*}
\La_{1,a} &= Y_{1,a}, \\
\La_{i,a} &= \La_{i-1,a} A_{i-1,aq^{i-1}}^{-1}, \quad \quad
i=2,\ldots,\el.
\end{align*}

$$
\begin{CD}
\bullet @>{1,q}>> \bullet @>{2,q^2}>> \bullet \quad \cdots \quad \bullet
@>{i,q^i}>> \bullet \quad \cdots \quad \bullet @>{\el,q^\el}>> \bullet
\end{CD}
$$

\medskip

\subsubsection{The $B_\el$ series}

$J = \{ 1,\ldots,\el,0,\ol{\el},\ldots,\ol{1} \}$.

\begin{align*}
\La_{i,a} &= Y_{i,aq^{2i-2}} Y_{i-1,aq^{2i}}^{-1}, \quad \quad
i=1,\ldots,\el-1, \\ \La_{\el,a} &= Y_{\el,aq^{2\el-3}}
Y_{\el,aq^{2\el-1}} Y_{\el-1,aq^{2\el}}^{-1}, \\ \La_{0,a} &=
Y_{\el,aq^{2\el-3}} Y_{\el,aq^{2\el+1}}^{-1}, \\ \La_{\ol{\el},a} &=
Y_{\el-1,aq^{2\el-2}} Y_{\el,aq^{2\el-1}}^{-1}
Y_{\el,aq^{2\el+1}}^{-1}, \\ \La_{\ol{i},a} &= Y_{i-1,aq^{4\el-2i-2}}
Y_{i,aq^{4\el-2i}}^{-1}, \quad \quad i=1,\ldots,\el-1.
\end{align*}

Equivalently,
\begin{align*}
\La_{1,a} &= Y_{1,a}, \\ \La_{i,a} &= \La_{i-1,a}
A_{i-1,aq^{2i-2}}^{-1}, \quad \quad i=2,\ldots,\el, \\ \La_{0,a} &=
\La_{\el,a} A_{\el,aq^{2\el}}^{-1}, \\ \La_{\ol{\el},a} &= \La_{0,a}
A_{\el,aq^{2\el-2}}^{-1}, \\ \La_{\ol{i},a} &= \La_{\ol{i+1},a}
A_{i,aq^{2(2\el-i-1)}}^{-1}, \quad \quad i=1,\ldots,\el-1.
\end{align*}

$$
\begin{CD}
\bullet @>{1,q^2}>> \bullet @>{2,q^4}>> \bullet \cdots \bullet
@>{\el,q^{2\el}}>> \bullet @>{\el,q^{2\el-2}}>> \bullet
@>{\el-1,q^{2\el}}>> \bullet \cdots \bullet @>{2,q^{4\el-6}}>> \bullet
@>{1,q^{4\el-4}}>> \bullet
\end{CD}
$$

\medskip

\subsubsection{The $C_\el$ series}

$J = \{ 1,\ldots,\el,\ol{\el},\ldots,\ol{1} \}$.

\begin{align*}
\La_{i,a} &= Y_{i,aq^{i-1}} Y_{i-1,aq^{i}}^{-1}, \quad \quad
i=1,\ldots,\el, \\ \La_{\ol{i},a} &= Y_{i-1,aq^{2\el-i+2}}
Y_{i,aq^{2\el-i+3}}^{-1}, \quad \quad i=1,\ldots,\el.
\end{align*}

Equivalently,
\begin{align*}
\La_{1,a} &= Y_{1,a}, \\
\La_{i,a} &= \La_{i-1,a} A_{i-1,aq^{i-1}}^{-1}, \quad \quad
i=2,\ldots,\el, \\
\La_{\ol{\el},a} &= \La_{\el,a} A_{\el,aq^{\el+1}}^{-1}, \\
\La_{\ol{i},a} &= \La_{\ol{i+1},a} A_{i,aq^{2\el-i+2}}^{-1},
\quad \quad i=1,\ldots,\el-1. 
\end{align*}

$$
\begin{CD}
\bullet @>{1,q}>> \bullet @>{2,q^2}>> \bullet \cdots \bullet
@>{\el-1,q^{\el-1}}>> \bullet @>{\el,q^{\el+1}}>> \bullet
@>{\el-1,q^{\el+3}}>> \bullet \cdots \bullet @>{2,q^{2\el}}>> \bullet
@>{1,q^{2\el+1}}>> \bullet
\end{CD}
$$

\medskip

\subsubsection{The $D_\el$ series}

$J = \{ 1,\ldots,\el,\ol{\el},\ldots,\ol{1} \}$.

\begin{align*}
\La_{i,a} &= Y_{i,aq^{i-1}} Y_{i-1,aq^{i}}^{-1}, \quad \quad
i=1,\ldots,\el-2, \\ \La_{\el-1,a} &= Y_{\el,aq^{\el-2}}
Y_{\el-1,aq^{\el-2}} Y_{\el-2,aq^{\el-1}}^{-1}, \\ \La_{\el,a} &=
Y_{\el,aq^{\el-2}} Y_{\el-1,aq^{\el}}^{-1}, \\ \La_{\ol{\el},a} &=
Y_{\el-1,aq^{\el-2}} Y_{\el,aq^\el}^{-1}, \\ \La_{\ol{\el-1},a} &=
Y_{\el-2,aq^{\el-1}} Y_{\el-1,aq^\el}^{-1} Y_{\el,aq^\el}, \\
\La_{\ol{i},a} &= Y_{i-1,aq^{2\el-i-2}} Y_{i,aq^{2\el-i-1}}^{-1},
\quad \quad i=1,\ldots,\el-2.
\end{align*}

Equivalently,
\begin{align*}
\La_{1,a} &= Y_{1,a}, \\ \La_{i,a} &= \La_{i-1,a}
A_{i-1,aq^{i-1}}^{-1}, \quad \quad i=2,\ldots,\el, \\
\La_{\ol{\el},a} &= \La_{\el-1,a} A_{\el,aq^{\el-1}}^{-1}, \\
\La_{\ol{\el-1},a} &= \La_{\ol{\el},a} A_{\el-1,aq^{\el-1}}^{-1} \\
&= \La_{\el,a} A_{\el,aq^{\el-1}}^{-1}, \\ \La_{\ol{i},a} &=
\La_{\ol{i+1},a} A_{i,aq^{2\el-i-2}}^{-1}, \quad \quad
i=1,\ldots,\el-2.
\end{align*}

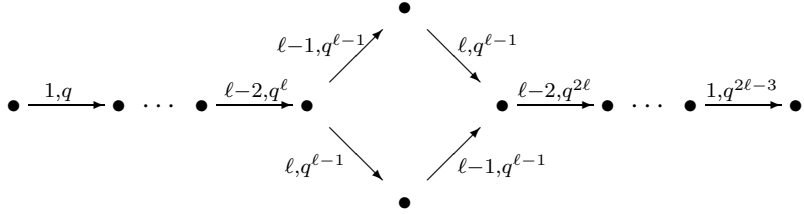
\setlength{\unitlength}{1mm}
\begin{center}
\begin{picture}(120,60)(0,0)
\put(0,30){$\bullet$}
\put(3,31){\vector(1,0){10}}
\put(14,30){$\bullet$}
\put(18,30){$\cdots$}
\put(25,30){$\bullet$}
\put(28,31){\vector(1,0){10}}
\put(39,30){$\bullet$}
\put(43,34){\vector(1,1){7}}
\put(52,43){$\bullet$}
\put(56,41){\vector(1,-1){7}}
\put(65,30){$\bullet$}
\put(43,28){\vector(1,-1){7}}
\put(52,17){$\bullet$}
\put(56,21){\vector(1,1){7}}
\put(68,31){\vector(1,0){10}}
\put(79,30){$\bullet$}
\put(83,30){$\cdots$}
\put(90,30){$\bullet$}
\put(93,31){\vector(1,0){10}}
\put(104,30){$\bullet$}
\put(5,32){${\scs 1,q}$}
\put(29,32){${\scs \el-2,q^{\el}}$}
\put(36,38){${\scs \el-1,q^{\el-1}}$}
\put(60,38){${\scs \el,q^{\el-1}}$}
\put(37,22){${\scs \el,q^{\el-1}}$}
\put(60,22){${\scs \el-1,q^{\el-1}}$}
\put(68,32){${\scs \el-2,q^{2\el}}$}
\put(93,32){${\scs 1,q^{2\el-3}}$}
\end{picture}
\end{center}


\section{Connection with Bethe Ansatz}    \label{Bethe}

As we have observed earlier \cite{FR:crit,FR:simple}, the formulas for
the $q$--characters have the same structure as the formulas for the
spectra of transfer-matrices on finite-dimen\-sional representations
of $\uqg$ obtained by the analytic Bethe Ansatz method
\cite{Re1,Re2,BR,KS1}. In this section we explain this connection in
more detail.

\subsection{Bethe Ansatz formulas}

Let us recall the idea of the Bethe Ansatz method (see the original
works \cite{FT,STF}). The problem is to find the eigenvalues of the
transfer-matrices $t_V(z), V \in \on{Rep} \uqg$, defined by formula
\eqref{tv}, on a finite-dimensional representation $W$. This problem
arises naturally, particularly, when $W = U^{\otimes N}$, in the study
of quantum spin chains, such as the XXZ model (see
\cite{FT,Re1,Re2}). Conjecturally, for each $V$, all eigenvalues of
$t_V(z)$ can be represented in a uniform way, and the eigenvalues are
parametrized by finite sets of complex numbers, which are solutions of
the so-called Bethe Ansatz equations (see \cite{Re1,Re2,BR,KS1} for
details).

We want to give an interpretation of these formulas from the point of
view of the $q$--characters. Suppose for simplicity that $W = V(\bp)$
is irreducible. Then it has a unique highest weight vector
$v_{\bp}$. Due to the property \eqref{hwv}, we obtain that $v_{\bp}$
is an eigenvector of $t_V(z)$, whose eigenvalue equals $\chi_q(V)$, in
which we substitute each $Y_{i,a}$ by $Y_i^{\bp}(za)$. We recall from
\secref{gen case} that $Y_i^{\bq}(za)$ can be found if we solve the
equations \eqref{express} for $i \in I$, in which we set
$$
A_{i,a} = q_i^{2+\deg P_i}
\dfrac{P_i(z^{-1}a^{-1}q_i^{-1})}{P_i(z^{-1}a^{-1}q_i)}.
$$
The Analytic Bethe Ansatz hypothesis is a statement that all other
eigenvalues of $t_V(z)$ on $\otimes_{j=1}^N V(\bp_j)$ have similar
structure.

\medskip

\noindent{\bf Hypothesis.} {\em Each eigenvalue of $t_V(z)$ on the
tensor product of finite-dimensional representations $\otimes_{i=1}^N
V(\bp_i)$ can be written as $\chi_q(V)$, in which we substitute
$Y_{i,a}$ by
\begin{equation}    \label{hyp}
\prod_{j=1}^N Y_i^{\bp_i}(za) \prod_{k=1}^{m_i}
\frac{za - w^{(i)}_k q_i}{za - w^{(i)}_k q_i^{-1}},
\end{equation}
for some $w^{(i)}_1,\ldots,w^{(i)}_{m_i}$. Furthermore, such an
eigenvalue occurs if an only if all seeming poles at $z=w^{(i)}_k
a^{-1} q_i^{-1}$ cancel each other.}

\medskip

The last condition can be written as a system of algebraic equations
on the zeroes $w_j$'s, which are called the Bethe Ansatz equations. It
turns out that the pole cancellation condition places such a stringent
constraint on the polynomial $\chi_q(V)$ that together with other
natural requirements, it suffices to reconstruct $\chi_q(V)$ completely
in many examples, see \cite{Re1,Re2,BR,KS1}.

Our construction of $q$--characters so far explains the Bethe Ansatz
hypothesis only for the vector $v_{\bp}$. However, such an explanation
can be given along the lines of our construction \cite{FFR} of the
spectra of the hamiltonians of the Gaudin system, which are certain
limits of the transfer-matrices $t_V(z)$ as $q \arr 1$. Recall that in
\cite{FFR} we constructed the eigenvectors of the Gaudin hamiltonians
using the Wakimoto modules at the critical level over the affine
Kac-Moody algebra $\G$. A generalization of this construction, in
which $\G$ is replaced by $\uqg$, should give us eigenvectors of
$t_V(z)$ with eigenvalues of the above form. The eigenvectors can be
constructed from intertwining operators between $\uqg$--modules at the
critical level (see Prop. 2.3 of \cite{DE}). To construct such
intertwiners explicitly, one has to use special singular vectors in
Wakimoto modules, which should exist precisely when the Bethe Ansatz
equations are satisfied, just as in the case of the Gaudin model
\cite{FFR}. We plan to discuss this in more detail in a separate
publication.

\subsection{The case of $\uqsl$}
According to formula \eqref{APhi}, we have:
$$
Y_{i,aq} Y_{i,aq^{-1}} = A_{i,a} = q^2 \Phi^-(z^{-1}a^{-1})
$$
(this is the equation \eqref{express} in this case). Note that the
value of $a$ in our formulas is not important, since we can always set
it to be $1$ by applying the automorphism $z \mapsto z a^{-1}$. Hence
from now on we will set $a=1$.

Let $v^{(r)}_0$ be the highest weight vector of $W_r(b)$. According to
formula \eqref{spec1}, we have
\begin{equation}
\Phi^-(u) \cdot v^{(r)}_0 = q^{r}
\frac{1-ubq^{-r}}{1-ubq^{r}} v^{(r)}_0.
\end{equation}
Therefore $Y^{(r)}(z) = Y^{P^{(r)}}(z)$ satisfies:
$$
Y^{(r)}(zq) Y^{(r)}(zq^{-1}) = q^{2+r}
\frac{1-z^{-1}bq^{-r}}{1-z^{-1}bq^{r}} = q^{2-r}
\frac{1-zb^{-1}q^{r}}{1-zb^{-1}q^{-r}},
$$
as a formal power series in $z$. This can be solved as follows:
\begin{equation}    \label{value}
Y^{(r)}(z) = q^{1-\frac{r}{2}} \frac{(zb^{-1}q^{r+1};q^4)_\infty
(zb^{-1}q^{-r+3};q^4)_\infty}{(zb^{-1}q^{r+3};q^4)_\infty
(zb^{-1}q^{-r+1};q^4)_\infty} =
\frac{\mu^{(r)}(zq^{-1})}{\mu^{(r)}(zq)},
\end{equation}
where
$$
\mu^{(r)}(z) =
\frac{(zb^{-1}q^{r+2};q^4)_\infty}{(zb^{-1}q^{-r+2};
q^4)_\infty}.
$$

The hypothesis above means in this case that the eigenvalues of
$t_{W_1(1)}$ on $\otimes_{j=1}^N W_{r_j}(b_j)$ are given by
\begin{equation}    \label{baxter}
q^\Delta \frac{\mu(zq^{-1})}{\mu(zq)} \frac{Q(zq^{-1})}{Q(zq)} +
q^{-\Delta} \frac{\mu(zq^3)}{\mu(zq)} \frac{Q(zq^3)}{Q(zq)},
\end{equation}
where
\begin{align*}
\mu(z) &= \prod_{j=1}^N
\frac{(zb_j^{-1}q^{r_j+2};q^4)_\infty}{(zb_j^{-1}q^{-r_j+2};
q^4)_\infty} = \prod_{j=1}^N \mu^{(r_j)}(zb_j^{-1}),
\\ Q(z) &= \prod_{k=1}^m (1-zw_k^{-1}),
\end{align*}
and $\Delta = 2m + \sum_j (1-r_j/2)$. Here $w_k$'s are complex
numbers subject to the equations:
$$
q^\Delta \frac{\mu(w_kq^{-2})}{\mu(w_k)} \prod_{s\neq k}
\frac{1-w_kw_s^{-1}q^{-2}}{1-w_kw_s^{-1}} + q^{-\Delta}
\frac{\mu(w_kq^2)}{\mu(w_k)} \prod_{s\neq k}
\frac{1-w_kw_s^{-1}q^{2}}{1-w_kw_s^{-1}} = 0,
$$
which mean that the expression \eqref{baxter} has no singularities at
the points $z=w_k q^{-1}$.

These equations are equivalent to the equations
\begin{equation}    \label{bae sl2}
\prod_{j=1}^N q^{r_j} \frac{w_k-b_j q^{-r_j}}{w_k-b_j q^{r_j}}
= - q^{-N} \prod_{s\neq k}
\frac{w_k-w_sq^{-2}}{w_k-w_s q^{2}}.
\end{equation}
These are the Bethe Ansatz equations for $\g=\sw_2$, which ensure that
the eigenvalues \eqref{baxter} have no poles at $z=w_k q^{-1}$.

Formula \eqref{baxter} first appeared in the works of R.~Baxter (see
\cite{Ba}), and polynomials $Q(z)$ are often called Baxter's
polynomials.

Note that $Q(z)=1-zw^{-1}$ is actually $\mu^{(r)}(z)$ in the special
case when $r=-2$. This reflects the special role that the Wakimoto
modules with highest weights $-\al_i$ play in the construction of the
eigenvectors (cf. \cite{FFR}).

\subsection{General case}

For a general Lie algebra $\g$ we find the Bethe Ansatz equations by
looking at the possible cancellations of poles at $z=w^{(i)} a^{-1}
q_i^{-1}$ in the sum of the expressions \eqref{hyp}. We find that
there are indeed possible pairwise cancellations between the
expressions \eqref{hyp} corresponding to monomials of the form $M$ and
$M A_{i,aq_i}^{-1}$. Each such cancellation gives rise to an equation,
which says that the sum of the residues of the expressions \eqref{hyp}
corresponding to $M$ and $M A_{i,aq_i}^{-1}$ at $z=w^{(i)} a^{-1}
q_i^{-1}$ is equal to $0$. Explicit calculation similar to the one
presented in the previous section gives us the following system of
equations (as before, we set $W = \bigotimes_{j=1}^N V(\bp_j)$, where
$\bp_j=(P_{j,1},\ldots,P_{j,\ell})$) on the variables $w^{(i)}_k$,
where $i=1,\ldots,\ell, k=1,\ldots,m_i$:
\begin{equation}    \label{bae gen}
\prod_{j=1}^N q_i^{\deg P_{j,i}}
  \frac{P_{j,i}(q_i^{-1}/w^{(i)}_k)}{P_{j,i}(q_i/w^{(i)}_k)} =
  - q^N \prod_{s \neq k}
  \frac{w^{(i)}_k-w^{(i)}_sq_i^{-2}}{w^{(i)}_k-w^{(i)}_s q_i^{2}}
  \prod_{l \neq i} \prod_{s=1}^{m_l}
  \frac{w^{(i)}_k-w^{(l)}_sq^{-C_{li}}}{w^{(i)}_k-w^{(l)}_s
  q^{C_{li}}}.
\end{equation}

This is the most general system of Bethe Ansatz equations
corresponding to an arbitrary collection of Drinfeld polynomials
$P_{j,i}, j=1,\ldots,N; i=1,\ldots,\ell$.

Because they come from ``local'' cancellations (those occurring
between monomials of the form $M$ and $M A_{i,aq_i}^{-1}$), these
equations are the {\em same} for each choice of the ``auxiliary''
space $V$ (which gives rise to the transfer matrix $t_V(z)$ acting on
the ``physical'' space $W$). Conjecturally, each solution gives rise
to a common eigenvector of all transfer-matrices $t_V(z), V \in
\on{Rep} \uqg$. The eigenvalues are conjecturally given by the linear
combinations of the expressions \eqref{hyp}, corresponding to the
$q$-characters $\chi_q(V)$, as stated in the above
Hypothesis.\footnote{Construction of the corresponding eigenvectors is
a separate question that will not be addressed here.} Moreover, we
expect that in a generic situation this Bethe Ansatz is ``complete'';
that is, the set of all solutions of the system \eqref{bae gen}
(modulo the obvious action of the product of symmetric groups) is in
one-to-one correspondence with the set of common eigevectors (or
eigenvalues, since we expect the spectrum to be simple, at least, in
the generic situation) of the transfer-matrices $t_V(z)$ on $W =
\bigotimes_{j=1}^N V(\bp_j)$.

For $\g=\sw_2$ and $V(\bp_j) = W_{r_j}(b)$, we have
$$
q^{\deg P_{j,1}}
  \frac{P_{j,1}(uq^{-1})}{P_{j,i}(uq)} = q^{r_j}
  \frac{1-ubq^{-r_j}}{1-ubq^{r_j}},
$$
and formula \eqref{bae gen} gives us formula \eqref{bae sl2}.

\section{The screening operators}    \label{screen}

In this section we define certain operators on $\yy =
\Z[Y_{i,a_i}^{\pm 1}]_{i\in I; a_i \in \C^\times}$, which we call the
screening operators. The terminology is explained by the fact that
these operators are certain limits of the screening operators used in
the definition of the deformed $\W$--algebras (see the next section).

\subsection{Definition of the screening operators}    \label{defscr}

Consider free $\yy$--module, which is the direct sum of vector spaces
$$
\wt{\yy}_i = \us{x \in \C^\times}{\opl} \yy \otimes S_{i,x}.
$$
Let $\yy_i$ be the quotient of $\wt{\yy}_i$ by the submodule generated
by elements of the form
\begin{equation}    \label{relation}
S_{i,xq_i^2} = A_{i,xq_i} S_{i,x}.
\end{equation}
Clearly,
$$
\yy_i \simeq \us{x \in (\C^\times/q_i^{2\Z})}{\opl} \yy \otimes S_{i,x},
$$
and so $\yy_i$ is also a free $\yy$--module.

Now define a linear operator $\wt{S}_i: \yy \arr \wt{\yy}_i$ by the
formula
$$
\wt{S}_i \cdot Y_{j,a} = \delta_{i,j} Y_{i,a} S_{i,a}
$$
and the Leibniz rule: $\wt{S}_i \cdot (ab) = (\wt{S}_i \cdot a) b + a
(\wt{S}_i \cdot b)$. In particular,
$$
\wt{S}_i \cdot Y_{j,a}^{-1} = - \delta_{i,j} Y_{i,a}^{-1} S_{i,a}.
$$

Finally, let $S_i: \yy \arr \yy_i$ be the composition of $\wt{S}_i$
and the projection $\wt{\yy}_i \arr \yy_i$. We call $S_i$ the $i$th
{\em screening operator}.

\begin{conj}    \label{main}
{\em The image of the homomorphism $\chi_q$ equals the intersection of
the kernels of the operators $S_i, i \in I$.}
\end{conj}

One can check explicitly that for all $\g$ of classical types, the
$q$--character of $V_{\omega_1(a)}$ lies in the intersection of the
kernels of the operators $S_i, i \in I$. In fact, in \cite{FR:simple}
that has been proved for the more general $(q,t)$--characters. This
means that the subring of $\on{Rep} \uqg$ generated by $t_1(a) =
\chi_q(V_{\omega_1(a)})$ lies in the intersection of the kernels of
the operators $S_i, i \in I$. If $\g = A_\el$ or $C_\el$, then this
subring coincides with $\on{Rep} \uqg$. For the $B_\el$ and $D_\el$
series, one needs to check that the $q$--characters of the spinor
representations also belong to the kernel of $S_i$'s. Explicit
formulas for the latter are given in \cite{FR:crit}.

For $x \in \C^\times$, denote by $\yy^{(x)}$ the subring
$\Z[Y_{i,xq^{2n_i}}]_{i \in I; n_i \in \Z}$ of $\yy$. Then we have:
$$
\yy \simeq \us{x \in (\C^\times/q^{\Z})}{\otimes} \yy^{(x)},
$$
and it follows from the definition of the operators $S_i$ that
$$
\bigcap_{i \in I} \on{Ker}_{\yy} S_i \simeq \us{x \in
(\C^\times/q^{\Z})}{\otimes} \bigcap_{i \in I} \on{Ker}_{\yy^{(x)}}
S_i.
$$
\conjref{main} implies that
$$
\on{Rep}^{(x)} \uqg \simeq \bigcap_{i \in I} \on{Ker}_{\yy^{(x)}} S_i,
$$
where $\on{Rep}^{(x)} \uqg = \Z[t_{i,xq^{n_i}}]_{i \in I; n_i \in
\Z}$.

Therefore we expect that any irreducible representation $V$ of $\uqg$,
considered as an element of $\Rep \uqg$, is the tensor product of
irreducible representations that belong to $\Rep^{(x)} \uqg$, so that
the structure of $\Rep \uqg$ is contained in the structure of
$\Rep^{(x)} \uqg$ (which are isomorphic to each other for different
$x$).

\subsection{Proof of \conjref{main} for $U_q \widehat{\sw}_2$}

In this case $\on{Rep} U_q \widehat{\sw}_2 \simeq \Z[t_a]_{a \in
\C^\times}$, where $t_a$ is the class of $V(a)$, which is the same as
$W_1(a)$ from \secref{casesl2}. The results of that section give us
the following formula:
$$
\chi_q(V) =  Y_a + Y_{aq^2}^{-1}.
$$
Thus, the image of $\chi_q$ is $\Z[Y_a + Y_{aq^2}^{-1}]_{a \in
\C^\times}$.

The action of the operator $S$ is given by the formula
$$
S \cdot Y_a^{\pm 1} = \pm Y_a^{\pm 1} S_a,
$$
and the Leibniz rule. We also have the following relation
$$
S_{aq^2} = S_a Y_a Y_{aq^2}.
$$

Now it is clear that
$$
\on{Ker}_{\yy} S = \us{x \in (\C^\times/q^{2\Z})}{\otimes}
\on{Ker}_{\yy^{(x)}} S,
$$
where $\yy^{(x)} = \Z[Y_{xq^{2n}}]_{n \in \Z}$.

On the other hand, it follows from \thmref{classl2} that
$$
\on{Rep} \uqsl = \us{x \in (\C^\times/q^{2\Z})}{\otimes}
\on{Rep}^{(x)} \uqsl,
$$
where $\on{Rep}^{(x)} \uqsl = \Z[t_{xq^{2n}}]_{n \in \Z}$.

Hence it suffices to show that
$$
\Z[t_{xq^{2n}}]_{n \in \Z} = \on{Ker}_{\yy^{(x)}} S.
$$

Without loss of generality we can set $x=1$. To simplify notation, we
denote $Y_{q^{2n}}$ by $y_n$ and $S_{q^{2n}}$ by $s_n$. We need to
prove that
$$
\Z[y_n + y_{n+1}^{-1}]_{n \in \Z} = \on{Ker} S,
$$
where
$$
S: \Z[y_n^{\pm 1}]_{n\in\Z} \arr \yy^{(1)} = \us{m \in \Z}{\opl}
\Z[y_n^{\pm 1}]_{n\in\Z} \otimes s_m/(s_m - y_m y_{m-1}
s_{m-1})
$$
is given by
$$
S \cdot y_n^{\pm 1} = \pm y_n^{\pm 1} s_n
$$
and the Leibniz rule. Note that $S$ commutes with the shift $y_n \arr
y_{n-k}$ for any integer $k$. Given an element in the kernel of $S$,
we can apply to it a shift with sufficiently large $k$ to make it into
an element of $\Z[y_n,y_{n+1}^{-1}]_{n\geq 0}$, which also belongs to
the kernel of $S$. Therefore without loss of generality we can
restrict ourselves to $\Z[y_n,y_{n+1}^{-1}]_{n\geq 0}$.

Furthermore, we can identify $\yy^{(1)}$ with $\Z[y_n^{\pm
1}]_{n\in\Z} \cdot s_0$ and hence with $\Z[y_n^{\pm
1}]_{n\in\Z}$. After this identification, $S$ becomes the derivation
$$
S \cdot y_n^{\pm 1} = \left\{ \begin{array}{cc} \pm y_n^{\pm 1}
\prod_{i=1}^n y_i y_{i-1}, \quad \quad n\geq 0 \\ \pm
y_n^{\pm 1} \prod_{i=0}^{n+1} y_i^{-1} y_{i-1}^{-1}, \quad \quad
n<0. \end{array} \right.
$$
In particular, we see that $\Z[y_n,y_{n+1}^{-1}]_{n\geq 0}$ is
$S$--invariant, and we want to show that the kernel of $S$ on
$\Z[y_n,y_{n+1}^{-1}]_{n\geq 0}$ equals $\Z[t_n]_{n\geq 0}$, where
$t_n = y_n + y_{n+1}^{-1}$.

Let us write
$$
\Z[y_n,y_{n+1}^{-1}]_{n\geq 0} = \Z[t_n,y_n]_{n\geq 0}/(t_n y_{n+1}
- y_n y_{n+1} - 1).
$$
Consider the set of monomials
$$
t_{n_1} \ldots t_{n_k} y_{m_1} \ldots y_{m_l},
$$
where all $n_i \geq 0$ and are lexicographically ordered, all $m_i
\geq 0$ are lexicographically ordered, and also $m_j \neq n_i+1$ for
all $i$ and $j$. We call these monomials reduced. It is easy to see
that the set of reduced monomials is a basis of
$\Z[y_n,y_{n+1}^{-1}]_{n\geq 0}$. Now let $P$ be an element of the
kernel of $S$. Let us write it as a linear combination of the reduced
monomials. We can then represent $P$ as $y_{N}^a Q + R$. Here $N$ is
the largest integer, such that $y_{N}$ is present in at least one of
the basis monomials appearing in its decomposition; $a$ is the largest
power of $y_{N}$ in $P$; $Q$ does not contain $y_{N}$, and $R$ is not
divisible by $y_{N}^a$; we assume here that both $y_{N} Q$ and $R$ are
linear combinations of reduced monomials.

When we apply $S$ to $P$ we obtain
\begin{equation}    \label{leading}
a y_{N}^{a+1} y_{N-1} \prod_{i=1}^{N-1} y_i y_{i-1} Q
\end{equation}
plus the sum of terms that are not divisible by $y_{N}^{a+1}$. Of
course, \eqref{leading} may not be in reduced form. But $Q$ does not
contain $t_{N-1}$. Therefore when we rewrite it as a linear
combination of reduced monomials, that linear combination will still
be divisible by $y_{N}^{a+1}$. On the other hand, no other terms in $S
\cdot P$ will be divisible by $y_{N}^{a+1}$. Hence for $P$ to be in
the kernel, \eqref{leading} has to vanish, which can only happen if
$P$ does not contain $y_m$'s at all, i.e., $P \in \Z[t_n]_{n\geq 0}$,
which is what we wanted to prove.

In the same way we obtain the following statement.

\begin{prop}    \label{keri}
The kernel of $S_i: \yy \arr \yy_i$ equals $${\mc R}_i =
\Z[Y_{j,a_j}]_{j\neq i; a_j \in \C^\times} \otimes \Z[Y_{i,b} +
Y_{i,b} A_{i,bq_i}^{-1}]_{b \in \C^\times}.$$
\end{prop}

\conjref{main} can therefore be interpreted as saying that the image
of $\chi_q$ in $\yy$ equals $\bigcap_{i \in I} {\mc R}_i$
(cf. \secref{gencase}).

\section{The connection with the deformed $\W$--algebras}
\label{defw}

Our motivation for the definition of the screening operators $S_i$ and
for \conjref{main} comes from the theory of deformed
$\W$--algebras. In this section we will explain this connection.

\subsection{The representation ring and the center at the critical
level}    \label{repcen}

We start by recalling the connection between $\on{Rep} \uqg$ and the
center of $\uqg$ at the critical level $-h^\vee$ (minus dual Coxeter
number).

Following \cite{RS,DE}, for each $V \in \on{Rep} \uqg$, we define in
addition to the $L$--operator $L_V(z)$ given by formula \eqref{lv},
the opposite $L$--operator
$$
L^-_V(z) = (\pi_{V(z)} \otimes \on{id})(\sigma(\R)),
$$
where $\sigma(a \otimes b) = b \otimes a$, and the total $L$--operator
$$
L^{\on{tot}}_V(z) = L_V(zq^{2h^\vee}) L^-_V(z).
$$
Let
$$
T_V(z) = \on{Tr}_V \; \left[ (q^{2\rho} \otimes 1) L^{\on{tot}}_V(z)
\right].
$$
Note that $T_V(z)$ is a formal power series, whose coefficients lie in
a completion of $\uqg$ and their action is well-defined on
any $\uqg$--module on which the action of the subalgebra $U_q \bb_+$
is locally finite. We will consider the projections of these elements
onto $\uqg_{\on{cr}} = \uqg/(c-q^{-h^\vee})$. Let $Z_q(\G)$ be the
center of $\uqg_{\on{cr}}$.

\begin{thm}[\cite{RS,DE}]
\hfill

{\em (1) For each $V \in \on{Rep} \uqg$, all Fourier coeficients of
$T_V(z)$ lie in $Z_q(\G)$.

(2) The map $V \arr T_V(z)$ is a $\C^\times$--equivariant ring
homomorphism $\wt{\nu}_q: \on{Rep} \uqg \arr Z_q(\G) \;\widehat{\otimes}\;
\C((z))$.}
\end{thm}

Now let $U_q \widehat{\n}_+$ be the subalgebra of $\uqg$ generated by
$x_i^+, i=0,\ldots,\el$. We will keep the same notation $U_q
\widehat{\n}_+$ and $U_q \bb_-$ for the projections of these
subalgebras onto $\uqg_{\on{cr}}$. Then we have the decomposition
$$
\uqg_{\on{cr}} = U_q \bb_- \otimes U_q \widehat{\n}_+,
$$
as a vector space, and so
$$
\uqg_{\on{cr}} = U_q \bb_- \oplus \left( U_q \bb_+ \otimes (U_q
\widehat{\n}_+)_0 \right),
$$
where $(U_q \widehat{\n}_+)_0$ is the augmentation ideal of $U_q
\widehat{\n}_+$. Denote by ${\mathbf p}$ be the corresponding
projection onto $U_q \bb_-$.

We find that
\begin{equation}    \label{multt}
{\mathbf p}(T_V(z)) = \on{Tr}_V \; q^{2\rho} L_V(z) T,
\end{equation}
where $T$ is defined by formula \eqref{t}. Thus, up to the inessential
factor of $T$, ${\mathbf p}(T_V(z))$ equals the transfer-matrix
$t_V(z)$ defined in formula \eqref{tv}.

Moreover, let ${\mathfrak z}_q(\G)$ be the subalgebra of $\uqb$
generated by the Fourier coefficients $t_V[n]$ of $t_V(z), V \in
\on{Rep} \uqg$. It is then easy to show that the restriction of
${\mathbf p}$ to $Z_q(\G)$ is a homomorphism of commutative algebras
$Z_q(\G) \arr {\mathfrak z}_q(\G)$, and we have the following
commutative diagram (up to the factor $T$ in \eqref{multt}):

\setlength{\unitlength}{1mm}
\begin{picture}(60,40)(-20,0)    \label{pic1}
\put(28,23){\vector(1,1){8}} 
\put(28,17){\vector(1,-1){8}} 
\put(47,29){\vector(0,-1){17}}
\put(11,19){$\on{Rep} \uqg$}
\put(37,32){$Z_q\left(\,\G\,\right)((z))$}
\put(37,6){${\mathfrak z}_{\,q}\left(\,\G\,\right)[[z]].$}
\put(48,20){${\mathbf p}$}    
\put(28,28){${\wt{\nu}_q}$}   
\put(32,14){${\nu_q}$}        
\end{picture}

\subsection{Free field realization of $\uqg$ and $q$--characters}

Recall that in \secref{hc} we defined the $q$--characters using the
Harish-Chandra homomorphism ${\mathbf h}_q$. A natural question is how
to construct an analogue of this homomorphism for $Z_q(\G)$. For that
one needs a free field (or Wakimoto) realization of $\uqg$. By this we
understand a family of homomorphisms ${\mc F}_k, k\in\C$, from $\uqg$
to the tensor product of a Heisenberg algebra ${\mc A}_q$ and the
Heisenberg algebra $\uqhh_{k+h^\vee}$, sending $c \in \uqg$ to
$q^k$. Here we denote by $\uqhh_{\al}$ the subalgebra of $\uqg$
generated by $c,k_i^\pm,h_{i,n}, i\in I, n\in \Z\backslash\{ 0 \}$,
modulo the relation $c=q^\al$.

Such a realization has been constructed in \cite{AOS} for $U_q
\widehat{\sw}_N$. Let us assume for a moment that it exists for any
$\uqg$.

At the critical level $k=-h^\vee$ the algebra $\uqhh_{0}$ is the
center of ${\mc A}_q \otimes \uqhh_{0}$, and therefore the image of
the center $Z_q(\G) \subset \uqg_{\on{cr}}$ under ${\mc F}_{-h^\vee}$
should lie in $\uqhh_{0}$. Furthermore, the commutative algebras
$Z_q(\G)$ and $\uqhh_{0}$ have natural Poisson structures, and the
resulting Harish-Chandra type homomorphism $\wt{\mathbf h}_q: Z_q(\G)
\arr \uqhh_{0}$ should preserve these Poisson structures. In the case
of $U_q \widehat{\sw}_N$, when the Wakimoto realization is available,
this homomorphism was analyzed in detail in \cite{FR:crit}, where it
was called the $q$--deformation of the Miura transformation.

Comparing the results of Sect.~3 with the results of \cite{FR:crit} in
the case of $U_q \widehat{\sw}_N$, we obtain the following commutative
diagram:
\begin{equation}    \label{pic2}
\begin{CD}
Z_q(\G)   @>{\wt{\mathbf h}_q}>>  \uqhh_{0}\\
@VV{\mathbf p}V      @VV{p}V\\
{\mathfrak z}_q(\G)   @>{{\mathbf h}_q}>>  \uqh
\end{CD}
\end{equation}
where $p: \uqhh_{0} \arr \uqh$ is the homomorphism that sends
$h_{i,n}, n>0$, to $0$, $h_{i,n}$ to $-h_{i,n}, n<0$, and $k_i$ to
$k_i^{-1}$. In particular, $p \circ \wt{\mathbf h}_q(T_V(z)) =
{\mathbf h}_q(t_V(z)) = \chi_q(V)$ (in which we replace $Y_{i,a}$ by
$Y_{i,a}^{-1}$).

We conjecture that the same is true in general. Thus, we expect that
the $q$--character homomorphism is a truncation of the free field
realization of the center $Z_q(\G)$ of $\uqg$ at the critical level.

The next step is to identify $Z_q(\G)$ with the classical limit of the
deformed $\W$--algebra.

\subsection{The algebra $\W_{q,t}(\g)$}    \label{wqt}

Let us recall the definition of the deformed
$\W$--algebra $\W_{q,t}(\g)$ and its free field realization from
\cite{FR:simple}.

Let $\HH(\g)$ be the Heisenberg algebra with generators $a_i[n],
i=1,\ldots,\el; n \in \Z$, and relations
\begin{equation}    \label{a}
[a_i[n],a_j[m]] = \frac{1}{n} (q^n - q^{-n}) (t^n - t^{-n}) \B_{ij}(q^n,t^n)
\delta_{n,-m},
\end{equation}
where
\begin{equation}    \label{qts}
\B_{ij}(q,t) = [\rr_i]_q \left( (q^{\rr_i} t^{-1} + q^{-\rr_i} t)
\delta_{i,j} - [I_{ij}]_q \right).
\end{equation}

There is a unique set of elements $y_i[n], i=1,\ldots,\el; n \in \Z$,
that satisfy:
\begin{equation}    \label{ay}
[a_i[n],y_j[m]] = \frac{1}{n} (q^{\rr_i n} - q^{-\rr_i n})(t^n - t^{-n})
\delta_{i,j} \delta_{n,-m}.
\end{equation}

Introduce the generating series:
$$
A_i(z) = q^{2a_i[0]}
:\exp \left( \sum_{m\neq 0} a_i[m] z^{-m} \right):,
$$
$$
Y_i(z) = q^{2 y_i[0]} :\exp \left( \sum_{m\neq 0} y_i[m] z^{-m}
\right):.
$$

Next we define the formal power series $S_i^+(z), i=1,\ldots,\el$, of
linear operators acting between Fock representations of $\HH(\g)$ that
satisfy the difference relations
\begin{equation}    \label{scr1}
S^+_i(zq^{-\rr_i}) = :A_i(z) S^+_i(zq^{\rr_i}):.
\end{equation}
These are the screening currents. Let $S^+_i$ be the $0$th Fourier
coefficient of $S^+_i(z)$.

\begin{rem} There is another set of screening currents: $S^-_i(z),
i=1,\ldots,\el$, introduced in \cite{FR:simple}, but we will not need
these currents here.\qed
\end{rem}

Let ${\mathbf H}_{q,t}(\g)$ be the vector space spanned by formal power
series of the form
\begin{equation}    \label{mono}
:\pa_z^{n_1} Y_{i_1}(zq^{j_1}t^{k_1})^{\ep_1} \ldots \pa_z^{n_1}
Y_{i_l}(zq^{j_l}t^{k_l})^{\ep_l}:,
\end{equation}
where $\ep_i = \pm 1$. The pair $(\hh,\pi_0)$, where $\pi_0$ is the
Fock representation of $\HH(\g)$ with highest weight $0$ (see
\cite{FR:simple}) is a deformed chiral algebra (DCA) in the sense of
\cite{FR:dca,FR:simple}.

We defined in \cite{FR:simple} the DCA ${\mathbf W}_{q,t}(\g)$ as the
maximal subalgebra of $(\hh,\pi_0)$, which commutes with the operators
$S^+_i, i=1,\ldots,\el$, i.e., the subspace of $\hh$, which consists
of all fields that commute with these operators. We also defined in
\cite{FR:simple} the deformed $\W$--algebra $\W_{q,t}(\g)$ as the
associative algebra, topologically generated by the Fourier
coefficients of fields from ${\mathbf W}_{q,t}(\g)$ (in the case of
$\sw_N$, the deformed $\W$--algebra had been previously constructed in
\cite{SKAO,FF:w,AKOS}, see also \cite{FR:crit,LP1}). All elements of
the algebra $\W_{q,t}(\g)$ act on the Fock representations $\pi_\la$
and commute with the screening operators.

\subsection{The classical limit of $\W_{q,t}(\g)$}

Now let us consider the limit of the algebras $\HH(\g), \W_{q,t}(\g)$
and the operators $S^+_i$ as $t \arr 1$. The algebra $\HH(\g)$ becomes
commutative with the Poisson structure given by
$$
\{ A,B \} = \lim_{t\arr 1} \frac{1}{2\log t} [A,B]_t.
$$
Using this formula, we obtain the following Poisson brackets between
the generators:
\begin{equation}    \label{poisson}
\{ a_i[n],a_j[m] \} = (q^{B_{ij} n} - q^{-B_{ij} n}) \delta_{n,-m},
\end{equation}
This formula shows that the map sending $a_i[n]$ to $-h_{i,n}$ is an
isomorphism $\hhh(\g) \simeq \uqhh_0$ of Poisson algebras.

The algebra $\W_{q,1}(\g)$ is the Poisson subalgebra of $\hhh(\g)$ that
consists of the elements Poisson commuting with the operators $S^+_i,
i=1,\ldots,\el$. Let us recall the following conjecture from
\cite{FR:simple}, which has been proved in the case of $U_q
\widehat{\sw}_N$ in \cite{FR:crit}.

\begin{conj}    \label{center}
{\em The Poisson algebra $\W_{q,1}(\g)$ is isomorphic to the Poisson
algebra $Z_q(\G)$. The embedding $\W_{q,1}(\g) \arr \hhh(\g)
\simeq \uqhh_0$ coincides with the free field realization homomorphism
of $Z_q(\G)$.}
\end{conj}

An analogous statement is true in the $q=1$ case as discussed in the
next section.

Now let us compute the Poisson brackets with $S^+_i$. Formula
\eqref{ay} becomes in the limit $t \arr 1$:
$$
\{ a_i[n],y_j[m] \} = (q^{\rr_i n} - q^{-\rr_i n}) \delta_{i,j}
\delta_{n,-m},
$$
which can be rewritten as
\begin{equation}    \label{genfun}
\{ A_i(z),Y_j(w)^{\pm 1} \} = \pm \left( \delta \left( \frac{wq_i}{z}
\right) - \delta \left( \frac{w}{zq_i} \right) \right) A_i(z)
Y_j(w)^{\pm 1}
\delta_{i,j},
\end{equation}
where
$$
\delta (x) = \sum_{n\in\Z} x^n.
$$

Formula \eqref{scr1} now reads:
\begin{equation}    \label{scrclass}
S^+_i(zq^{-\rr_i}) = A_i(z) S^+_i(zq^{\rr_i}).
\end{equation}
Combining it with \eqref{genfun} we obtain:
$$
\{ S^+_i(z),Y_j(w)^{\pm 1} \} = \pm \delta \left( \frac{w}{z} \right)
S^+_i(z) Y_j(w)^{\pm 1} \delta_{i,j} = \pm \delta \left( \frac{w}{z}
\right) Y_j(w)^{\pm 1} S^+_i(w) \delta_{i,j}.
$$
Taking the $0$th Fourier coefficient in $z$ we obtain:
\begin{equation}    \label{scraction}
\{ S^+_i,Y_j(w)^{\pm 1} \} = \pm Y_j(w)^{\pm 1} S_i(w) \delta_{i,j}.
\end{equation}

Now let $\yyy = \C[Y_j(wq^{2n})^{\pm 1}]_{j=1,\ldots,\el;n \in
\Z}$ and
\begin{align*}
\yyy_i &= \left( \us{m\in\Z}{\opl} \C[Y_j(wq^{2n})^{\pm 1}] \otimes
S_i(wq_i^{2m}) \right)/(S^+_i(zq_i^{-1}) - A_i(z) S^+_i(zq_i)) \\ &
\simeq \C[Y_j(wq^{2n})^{\pm 1}] \otimes S_i(w)
\end{align*}
Then $\{ S^+_i,\cdot \}$ is a linear operator $\yyy \arr
\yyy_i$.

The Poisson algebra $\hhh(\g)$ consists of the Fourier coefficients of
elements of $\yyy$. Let $\hhh(\g)_i$ be the span of the Fourier
coefficients of elements of $\yyy_i$. Then $\{ S^+_i,\cdot \}$
gives rise to an operator $\hhh(\g) \arr \hhh(\g)_i$. \conjref{center}
then tells us that
$$
Z_q(\G) \simeq \bigcap_{i=1,\ldots,\el} \on{Ker}_{\hhh(\g)} \{
S^+_i,\cdot \}.
$$

On the other hand, we can consider $\yyy$ itself as the classical
limit of the deformed chiral algebra ${\mathbf H}_{q,t}(\g)$ -- the
``space of fields'' of $\hhh(\g)$, and its subspace
$$\bigcap_{i=1,\ldots,\el} \on{Ker}_{\yyy} \{ S^+_i,\cdot \}$$ as
the classical limit of ${\mathbf W}_{q,t}(\g)$ -- the ``space of
fields'' of the $\W$--algebra $\W_{q,1}(\g)$. Then \conjref{center}
and the commutative diagrams \eqref{pic1}, \eqref{pic2} suggest the
following isomorphism
\begin{equation}    \label{main1}
\bigcap_{i=1,\ldots,\el} \on{Ker}_{\yyy} \{ S^+_i,\cdot \}
\simeq \on{Rep}^{(w)} \uqg.
\end{equation}

If we identify $Y_j(wq^{2n})$ with $Y_{j,q^{2n}}^{-1}$, then $\yyy$
and $\yyy_i$ gets identified with $\yy$ and $\yy_i$, respectively, and
the operator $\{ S^+_i,\cdot \}$ becomes the operator $S_i$ introduced
in \secref{defscr}. Hence \eqref{main1} is equivalent to
\conjref{main}.

\subsection{Comparison with the case $q=1$}

Now we explain the analogous picture in the case of affine Lie
algebras, i.e., when $q=1$. In this case most of the conjectures of
the previous subsections become theorems, except that for $q=1$ there
is no obvious connection between the center $Z(\G)$ and $\on{Rep} \G$.

\subsubsection{The definition of conformal $\W$--algebras.}

Let $\Hh(\g)$ be the Heisenberg algebra with generators $\ab_i[n],
i=1,\ldots,\el; n \in \Z$, and relations:
\begin{equation}    \label{oh}
[\ab_i[n],\ab_j[m]] = n \B_{ij} \beta \delta_{n,-m}.
\end{equation}
The conformal screening currents satisfy the differential equations
\begin{equation}    \label{diff1}
\pa_z \scr^+_i(z) = - \frac{1}{r_i} :\Ab_i(z) \scr^+_i(z):,
\end{equation}
where
$$
\Ab_i(z) = \sum_{n\in\Z} \ab_i[n] z^{-n-1}.
$$
The conformal screening operators are given by
$$
\scr_i^\pm = \int \scr_i^\pm(z) dz.
$$

Let $\pi_0$ be the vertex operator algebra (VOA) associated to the
Heisenberg algebra $\Hh(\g)$ (see \cite{FF:laws}). For generic $\beta$
the VOA ${\mathbf W}_\beta(\g)$ is defined \cite{FF:ds,FF:laws} as the
vertex operator subalgebra of the VOA $\pi_0$, which is the
intersection of kernels of the screening operators $\scr_i^-,
i=1,\ldots\el$:
$$
{\mathbf W}_\beta(\g) = \bigcap_{i=1,\ldots,\el} \on{Ker}_{\pi_0}
\scr_i^+.
$$
Thus, ${\mathbf W}_\beta(\g)$ consists of the fields that commute with
$\scr_i^+$. The $\W$--algebra $\WW(\g)$ is defined as the associative
(or Lie) algebra generated by the Fourier coefficients of these
fields.

The $\W$--algebra $\W_\beta(\g)$ can be obtained from $\W_{q,t}(\g)$
as $q \arr 1$ with $t = q^\beta$ (see \cite{FR:simple}).

\subsubsection{Classical limit}
The classical limit of $\W_\beta(\g)$ as $\beta \arr 0$ coincides with
a limit of $\W_{q,1}(\g)$ as $q \arr 1$. In this limit, the algebra
$\Hh(\g)$ becomes commutative and it inherits a Poisson structure,
which is given on the generators by the formula
\begin{equation}    \label{ohclass}
\{ \ab_i[n],\ab_j[m] \} = n \B_{ij} \delta_{n,-m}.
\end{equation}

Introduce new variables $\yb_i[n]$, such that
$$
\{ \ab_i[n],\yb_j[m] \} = n \delta_{i,j} \delta_{n,-m}.
$$
We rewrite this as
\begin{equation}    \label{AY}
\{ \Ab_i(z),\Yb_j(w) \} = - \delta_{i,j} w^{-1} \pa_z \delta \left(
\frac{w}{z} \right).
\end{equation}

In the limit $\beta \arr 0$, formula \eqref{diff1} becomes
\begin{equation}    \label{classdiff1}
\pa_z \scr^+_i(z) = - \frac{1}{r_i} \Ab_i(z) \scr^+_i(z)
\end{equation}
Thus, we can consider $\scr^+_i(z)$ as $\exp \left( - \Phi_i^\vee(z)
\right)$, where $\pa_z \Phi_i^\vee(z) = \Ab_i(z)/r_i$.

Combining \eqref{AY} and \eqref{classdiff1} we obtain:
$$
\{ \scr^+_i(z),\Yb_j(w) \} = \frac{1}{r_i} \delta_{i,j} w^{-1}
\delta \left( \frac{w}{z} \right) \scr^+_i(w),
$$
and hence
\begin{equation}    \label{dual}
\{ \scr^+_i,\pa_w^n \Yb_j(w) \} = \frac{1}{r_i} \pa_w^n \scr^+_i(w)
\delta_{i,j},
\end{equation}
where $\pa_w^n \scr^+_i(w)$ can be rewritten as a differential
polynomials in $\Ab_i(w)$ times $\scr^+_i(w)$ using formula
\eqref{classdiff1}.

Let $\pi_0 = \C[\ab_i[n]]_{n\leq -1}$ be the Fock representation of
${\mc H}_0(\g)$, on which $\ab_i[n]$ acts as multiplication by itself
for $n\leq -1$ and by $0$ for $n\geq 0$. We can identify $\pi_0$
naturally with ${\mc U} = \C[\pa_z^m \Ab_i(z)]_{i=1\,\ldots,\el; m\geq
0}$ by sending $P(\pa_z^m \Ab_i(z))$ to $P(\pa_z^m \Ab_i(z)) \cdot
1|_{z=0} \in \pi_0$. This identification, which is the remnant of the
VOA structure on $\pi_0$ for $\beta \neq 0$ is actually a ring
isomorphism sending $\ab_i[n], n\leq -1$, to $\pa_z^{-n-1}
\Ab_i(z)/(-n-1)!$.

Furthermore, the $\beta$--linear term of the VOA structure defines on
$\pi_0$ the structure of coisson algebra \cite{BD} (also called vertex
Poisson algebra in \cite{EF}) and the classical limit of the VOA
${\mathbf W}_\beta(\g)$, ${\mathbf W}_0(\g) \subset \pi_0$, can be
viewed as a coisson subalgebra of $\pi_0$. The coisson structure
encodes the Poisson structures in ${\mc H}_0(\g)$ and $\W_0(\g)$,
which are spanned by the Fourier coefficients of fields from $\pi_0$
and ${\mathbf W}_0(\g)$, respectively.

Let now ${\mc U}_i$ be the free ${\mc U}$--module with generator
$S^+_i(w)$. We extend the action of $\pa_w$ to ${\mc U}_i$ using
formula \eqref{classdiff1}. Define the operator $\scr^+_i: {\mc U}
\arr {\mc U}_i$ by formulas \eqref{dual}, \eqref{classdiff1}, the
Leibniz rule and the property that it commutes with the action of
$\pa_w$ on both spaces. Then by definition
$$
{\mathbf W}_0(\g) = \bigcap_{i=1,\ldots,\el} \on{Ker}_{\mc U} \{
\scr^+_i,\cdot \}.
$$

Let ${\mc H}_0(\g)_i$ be the span of the Fourier coefficients of
fields from ${\mc U}_i$. The linear operator $\{ \scr^+_i,\cdot \}$
acts from ${\mc H}_0(\g)$ to ${\mc H}_0(\g)_i$, and by definition the
classical $\W$--algebra is
\begin{equation}    \label{wo}
\W_0(\g) = \bigcap_{i=1,\ldots,\el} \on{Ker}_{{\mc H}_0(\g)} \{
\scr^+_i,\cdot \}.
\end{equation}

\begin{rem}
What we denote by $\WW(\g)$ here is really $\W_{-\sqrt{\beta}}(\gL)$ in
the notation of \cite{FF:laws}, where $\gL$ stands for the Langlands
dual Lie algebra to $\g$. The reason for the appearance of $\gL$ is
the factor of $1/r_i = 2r^\vee/(\al_i,\al_i)$ in formula
\eqref{diff1}. While generators $\ab_i[n]$ correspond to the simple
roots of $\g$, the rescaled generators $\ab_i^\vee[n] = \ab_i[n]/r_i$
correspond to the coroots of $\g$ and hence to the roots of $\gL$.

In particular, what we denote by $\W_0(\g)$ here is really $\W(\gL)$
of \cite{FF:laws}, which is the Poisson algebra of integrals of motion
of the Toda field theory associated to $\gL$. Again, this is due to
the factor $1/r_i$ in formula \eqref{classdiff1}, which makes the sum
$$
\sum_{i=1}^\el \scr^+_i = \sum_{i=1}^\el \int \exp (- \Phi_i^\vee(z))
dz
$$
the hamiltonian of the Toda field theory of $\gL$, see \cite{FF:laws}
for more details. This Poisson algebra can also be obtained by the
Drinfeld-Sokolov reduction \cite{DS} from the dual space to the Lie
algebra $\widehat{\gL}$, see \cite{FF:ds}. The embedding $\W_0(\g) \arr
{\mc H}_0(\g)$ is the classical Miura transformation.\qed
\end{rem}

\subsubsection{The center}

Now consider the affine Kac-Moody algebra $\G = \g[t,t^{-1}] \oplus \C
K$ and the completion $\wt{U}(\G)_k$ of $U(\G)/(K+k)$ introduced in
\cite{FF:ds}. Let $Z(\G)$ be the center of $\wt{U}(\G)_{\on{cr}} =
\wt{U}(\G)_{-h^\vee}$. We have the following result.

\begin{thm}[\cite{FF:ds}]
{\em $Z(\G)$ is isomorphic to $\W_0(\g)$ as a Poisson algebra.}
\end{thm}

This theorem and \eqref{wo} imply that
$$
Z(\G) \simeq \bigcap_{i=1,\ldots,\el}
\on{Ker}_{{\mc H}_0(\g)} \{ \scr^+_i,\cdot \}.
$$

The classical Miura transformation $\W_0(\g) \arr {\mc H}_0(\g)$ can be
interpreted in terms of the free field (Wakimoto) realization of
$\G$. Recall \cite{W,FF:wak} that this is a homomorphism from
$\wt{U}(\G)_k$ to the tensor product of the Heisenberg algebras ${\mc
A}$ and $\uhh_{k+h^\vee}$. The latter is the completion of the
universal enveloping algebra of the homogeneous Heisenberg subalgebra
$\widehat{\h} = \h[t,t^{-1}] \oplus \C K$ of $\G$ in which we set the
central element $K$ to be equal to $k+h^\vee$. In particular, when
$k=-h^\vee$, it becomes commutative, and inherists a Poisson
structure. The resulting Poisson algebra is isomorphic to ${\mc
H}_0(\g)$.

Under the free field homomorphism, the image of $Z(\G) \subset
\wt{U}(\G)_{\on{cr}}$ lies in $\uhh_0$. The resulting homomorphism
$\wt{{\mathbf h}}: Z(\G) \arr \uhh_0$ preserves the Poisson structures
and we have the following commutative diagram (see \cite{FF:ds}):
\begin{equation}    \label{pic3}
\begin{CD}
Z(\G)   @>{\wt{\mathbf h}}>>  \uhh_{0}\\
@VVV      @VVV\\
\W_0(\g)   @>{{\mathbf h}}>>  {\mc H}_0(\g)
\end{CD}
\end{equation}
where the isomorphism $\uhh_{0} \arr {\mc H}_0(\g)$ sends $h_{i,n}$ to
$-\ab_i[n]$. The bottom arrow is the Miura transformation $\W_0(\g)
\arr {\mc H}_0(\g)$, see \cite{FF:ds,FR:crit} for more detail.

Finally, there is an analogue of the commutative diagram \eqref{pic2}.
Let $\wt{\h} = \h \otimes t^{-1}\C[t^{-1}] \subset \h[t,t^{-1}] =
\widehat{\h}$, and $p$ be the quotient of $\uhh_0$ by the ideal
generated by $\h[t]$. Denote $\G_- = \g \otimes t^{-1}\C[t^{-1}]$. In
the same way as in the quantum case we define the projections
${\mathbf p}: \wt{U}(\G)_{\on{cr}} \arr U \G_-$ and ${\mathbf h}: U
\bb_- \arr \uh$. Let ${\mathfrak z}(\G)$ be the commutative subalgebra
of $U \g_-$ that is the image of $Z(\G)$ under ${\mathbf p}$. Then we
have:
\begin{equation}    \label{pic4}
\begin{CD}
Z(\G)   @>{\wt{\mathbf h}}>>  \uhh_{0}\\
@VV{\mathbf p}V      @VV{p}V\\
{\mathfrak z}(\G)   @>{{\mathbf h}}>>  \uh
\end{CD}
\end{equation}
up to the automorphism of $\uh$ that sends $h_{i,n}$ to $-h_{i,n}$.

Consider now the homomorphism of coisson algebras $\uh \arr \pi_0$
sending $h_{i,n}$ to $\ab_i[n]$. \thmref{center} can be rephrased as
saying that the corresponding map ${\mathfrak z}(\g) \arr \pi_0$ is an
isomorphism onto ${\mathbf W}_0(\g) \subset \pi_0$.

\subsubsection{Example. The case of $\widehat{\sw}_2$}

Let $\{ e,h,f \}$ be the standard basis of $\sw_2$. Set $a_n = a
\otimes t^n$ and
$$
a(z) = \sum_{n\in\Z} a_n z^{-n-1}.
$$
Introduce the generating function of the Sugawara operators $\wt{S}_n$
by formula
$$\wt{S}(z) = \sum_{n \in \Z} \wt{S}_n z^{-n-2} = \frac{1}{4} :h(z)^2:
+ \frac{1}{2} :e(z) f(z): + \frac{1}{2} :f(z) e(z):.$$ It is
well-known that $\wt{S}_n \in Z(\su)$, and $Z(\su)$ is topologically
generated by $\wt{S}_n, n \in \Z$.

The Lie subalgebra $\G_-$ of $\su$ is spanned by $e_n,h_n,f_n, n <
0$. We find:
$$
{\mathbf p}(\wt{S_n}) = \begin{cases} S_n, & n\leq -2 \\ 0, & n > -2,
\end{cases}
$$
where
$$
S_n = \frac{1}{4} \sum_{k+m=n; k,m<} h_m h_k + \frac{1}{2}
\sum_{k+m=n; k,m<0} \left( e_k f_m + f_m e_k \right).
$$
Hence we obtain:
\begin{equation}    \label{sn}
{\mathbf h}(S_n) = \frac{1}{4} \sum_{k+m=n; k,m\leq 0} h_m h_k -
\frac{1}{2} (n+1) h_n, \quad \quad n\leq -2.
\end{equation}
In fact, ${\mathfrak z}(\su) = \C[S_n]_{n\leq -2}$ and formula
\eqref{sn} defines its embedding into $\uh = \C[h_n]_{n\leq -1}$.

On the other hand, let $\chi_n, n\in\Z$, be the generators of
$\widehat{\h}$. Then one finds (see \cite{FR:crit}):
\begin{equation}    \label{tildesz}
\wt{{\mathbf h}}(\wt{S}(z)) = \frac{1}{4} \chi(z)^2 - \frac{1}{2}
\pa_z \chi(z).
\end{equation}

Formulas \eqref{sn} and \eqref{tildesz} show that the diagram
\eqref{pic4} is commutative up to the automorphism of $\uh$ sending
$h_n$ to $-h_n$.

Now consider the $\W$--algebra $\W_0(\sw_2)$, which is in this case
the classical Virasoro algebra. We have: ${\mc U} = \C[\pa_z^n
\Yb(z)]_{n\geq 0}$ and ${\mc U}_1 = {\mc U} \otimes S^+(z)$. The
operator $\{ \scr^+,\cdot \}: {\mc U} \arr {\mc U}_1$ is given by
the formula
$$
\{ \scr^+,\pa_z^n \Yb(z) \} = \pa_z^n \scr^+(z),
$$
the relation
$$
\pa_z \scr^+(z) = - 2 \Yb(z) \scr^+(z),
$$
and the Leibniz rule.

One finds \cite{FF:laws} that ${\mathbf W}_0(\sw_2) = \on{Ker} \{
\scr^+,\cdot \} = \C[\pa_z^n {\mathbf T}(z)]_{n\geq 0}$, where
\begin{equation}    \label{tz}
{\mathbf T}(z) = \Yb(z)^2 + \pa_z \Yb(z).
\end{equation}

Formulas \eqref{tildesz} and \eqref{tz} show the commutativity of
diagram \eqref{pic3} (note that $\Yb(z)$ goes to $-\frac{1}{2}
\chi(z)$).

\subsection{Deformed $\W$--algebra as a quantization of $\on{Rep} \uqg$}

According to \secref{wqt}, the DCA ${\mathbf W}_{q,t}(\g)$ can be
viewed as a quantization (along the $t$ variable) of the
representation ring $\on{Rep} \uqg$, while the $\W$--algebra
$\W_{q,t}(\g)$ is a quantization of the center $Z_q(\G)$.

The structure of the DCA ${\mathbf W}_{q,t}(\g)$ is similar to that of
the representation ring $\on{Rep} \uqg$. In particular, each
finite-dimensional representation $V$ of $\uqg$ should give rise to a
$(q,t)$--character $\chi_{q,t}(V) \in {\mathbf W}_{q,t}(\g)$, which
becomes $\chi_q(V)$ at $t=1$. These $(q,t)$--characters should be
linear combinations of normally ordered monomials in $Y_i(zq^mt^n)$
whose coefficients are rational functions in $q$ and $t$ taking
non-negative integral values at $t=1$ and $q=\ep,1$ (see
\cite{FR:simple}). The $(q,t)$--characters for the first fundamental
representation of $\uqg$, where $\g$ is of classical type are given in
\cite{FR:simple}; for fundamental representations of $U_q G^{(1)}_2$
they are given in \cite{BP} (see also \cite{Ko}).

Tensor product structure of $\on{Rep} \uqg$ is reflected in the pole
structure of the fusion of the $(q,t)$--characters. Namely, in all
known examples (see \cite{FR:dca,FR:simple,BP}), whenever $U$ appears
in the decomposition of the tensor product $V \otimes W$, there exists
an integer $a$, such that $\on{Res}_{z=wt^a} \chi_{q,t}(V)(z)
\chi_{q,t}(W)(w) dz/z$ equals $\chi_{q,t}(U)(w)$ up to a constant
multiple. It is natural to conjecture that this is true in general.

Now consider other classical limits of ${\mathbf W}_{q,t}(\g)$ and
$\W_{q,t}(\g)$ (see \cite{FR:simple}).

If $\g$ is simply-laced, then the other classical limit is $q \arr
1$. But $\W_{q,t}(\g)$ is actually invariant under the replacement $q
\arr t^{-1}, t \arr q^{-1}$ in this case, and so the structure of
$\W_{1,t}(\g)$ is essentially the same as that of $\W_{q,1}(\g)$.

If $\g$ is non-simply laced, then there are two interesting limits: $q
\arr 1$, and $q \arr \ep = \exp(\pi i/r^\vee)$.

According to the \cite{FR:simple}, the algebra ${\mathcal
W}_{\ep,t}(\g)$ is not commutative, but it contains a large center
${\mathcal W}'_t(\g)$. Conjecture 4 of \cite{FR:simple} states that
${\mathcal W}'_t(\g)$ is isomorphic to the center of $U_t(\GL)$ at the
critical level, where $\GL$ is the twisted affine algebra that is
Langlands dual to $\G$ (its Dynkin diagram is obtained from the Dynkin
diagram of $\G$ by reversing the arrows).

Alternatively, one can say that the limit $q \arr \ep$ of ${\mathbf
W}_{q,t}(\g)$ contains a commutative subalgebra that is isomorphic to
$\on{Rep} U_t(\GL)$. In that sense, ${\mathbf W}_{q,t}(\g)$ appears to
be a simultaneous quantization of $\on{Rep} \uqg$ and $\on{Rep}
U_t(\GL)$. Some evidence for this is presented in \cite{FR:simple},
where the $(q,t)$--characters $\chi_{q,t}(V_{\omega_1}(z))$ are given
explicitly for $\g$ of classical type. If we set $t=\ep$ in those
formulas, we indeed obtain the $t$--characters of representations of
$U_t(\GL)$, which can be seen by comparing them with the corresponding
formulas for the spectra of transfer-matrices \cite{Re2,KS2}. The
$(q,t)$--characters for the fundamental representations in the case
$\g=G_2$ obtained in \cite{BP} also agree with the conjecture.

Note that our conjecture above on the residues in the operator product
expansion of $(q,t)$--characters suggests that even the tensor
structures of $\on{Rep} \uqg$ and $\on{Rep} U_t(\GL)$ are related.

The analysis of the above formulas also shows that in the limit $q
\arr 1$ they look like characters of some algebra closely related to
$U_t(^L\; \widehat{\gL})$ (see \cite{FR:simple}). The limit $q \arr 1$
can be described alternatively using the difference Drinfeld-Sokolov
reduction, which we now recall.

\subsection{Difference Drinfeld-Sokolov reduction}

According to Conjecture 3 of \cite{FR:simple}, ${\mathbf W}_{1,t}(\g)$
and $\W_{1,t}(\g)$ can be obtained by the $t$--difference
Drinfeld-Sokolov reduction from the loop group of $G$
\cite{FRS,SS}. On the other hand, as we remarked above, if $\g$ is
simply-laced, then $\on{Rep} U_t(\G)$ equals ${\mathbf
W}_{1,t}(\g)$. Therefore for simply-laced $\g$, $\on{Rep} \uqg$ can be
obtained by the $q$--difference Drinfeld-Sokolov reduction. This
gives one an alternative method to find the $q$--characters of
irreducible representations with a given highest weight.

Consider for example the case of $\sw_N$ (see \cite{FRS}). Let
$M_{n,q}$ be the vector space of first order difference operators $D +
A(s)$, where $A(s)$ is an element of the formal loop group
$LSL_n=SL_n((s))$ of $SL_n$, and $D$ is the difference operator $(D
\cdot f)(s) = f(sq^2)$. The group $LSL_n$ acts on this manifold by the
$q$--gauge transformations \beq
\label{qadjoint} g(s) \cdot (D + A(s)) = g(sq^2)(D +
A(s))g(s)^{-1},
\end{equation}
i.e. $g(s) \cdot A(s) = g(sq^2) A(s) g(s)^{-1}$.

Now we consider the submanifold $M^J_{n,q} \subset M_{n,q}$ which
consists of operators $D+A(s)$, where $A(s)$ is of the form
\beq
\begin{pmatrix}
* & * & * & \hdots & * & * \\
-1 & * & * & \hdots & * & * \\
0 & -1 & * & \hdots & * & * \\
\hdotsfor{6} \\
0 & 0 & 0 & \hdots & * & * \\
0 & 0 & 0 & \hdots & -1 & *
\end{pmatrix}.
\end{equation}
It is preserved under the $q$--gauge action of the
group $LN$.

\begin{lem}    \label{qfree}
The action of $LN$ on $M^J_{n,q}$ is free and each orbit contains a
unique operator of the form \beq \label{qcan} D +
\begin{pmatrix} t_1 & t_2 & t_3 & \hdots & t_{n-1} & 1 \\ -1 & 0 & 0 &
\hdots & 0 & 0 \\ 0 & -1 & 0 & \hdots & 0 & 0 \\ \hdotsfor{6} \\ 0 & 0
& 0 & \hdots & 0 & 0 \\ 0 & 0 & 0 & \hdots & -1 & 0
\end{pmatrix}.
\end{equation}
\end{lem}

For other groups, the analogous statement has been proved by
Semenov-Tian-Shan\-sky and Sevostyanov \cite{SS}.

Now let $\Fq$ be the vector space the $q$--difference operators
\beq    \label{Lambda}
\La = D + \begin{pmatrix}
\la_1(s) & 0 & \hdots & 0 & 0 \\
-1 & \la_2(sq^{-2}) & \hdots & 0 & 0 \\
\hdotsfor{5} \\
0 & 0 & \hdots & \la_{n-1}(sq^{-2n+4}) & 0 \\
0 & 0 & \hdots & -1 & \la_n(sq^{-2n+2})
\end{pmatrix},
\end{equation}
where $\prod_{i=1}^n \la_i(sq^{-2i+2}) = 1$.

Let $\mu_q: \Fq \arr \Mq$ be the composition of the embedding
$\Fq \arr M^J_{n,q}$ and $\pi_q: M^J_{n,q} \arr M^J_{n,q}/LN \simeq
\Mq$. Using the definition of $\pi_q$ above, one easily finds that for
$\La$ given by \eqref{Lambda}, $\mu_q(\La)$ is the operator
\eqref{qcan}, where $t_i(s)$ is given by
\beq    \label{formulai}
t_i(s) = \sum_{j_1 < \ldots < j_i} \la_{j_1}(s) \la_{j_2}(sq^{-2})
\ldots \la_{j_i}(sq^{-2i+2}).
\end{equation}
This coincides with the formula for the $q$--character of the $i$th
fundamental representation of $U_q \sun$, if we make the replacement
$\la_i(s) \arr \La_{i,s}$, where $\La_{i,s}$ is defined in
\secref{exa}.  Thus, if we consider $\mu_q$ as a homomorphism
$$
\C[t_i(s)]_{i=1,\ldots,n-1;s\in \C^\times} \arr
\C[\la_i(s)]_{i=1,\ldots,n;s\in \C^\times}/\left( \prod_{i=1}^n
\la_i(sq^{-2i+2}) - 1 \right),
$$
then it coincides with the $q$-character homomorphism $\chi_q$. For
$q=1$, this result is due to Steinberg \cite{St}.

Now we see that the problem of reconstructing the $q$--character from
the dominant monomial $Y_{i_1,a_1} \ldots Y_{i_k,a_k}$ is equivalent
to finding the minimal combination of monomials in $\la_i(sq^{2n_i})$
with positive integral coefficients, which lies in the image of the
homomorphism $\mu_q$. One can probably use the geometry of the orbit
space $M^J_{n,q}/LN$ to study this question. This method can also be
applied to other simply-laced $\g$.

\end{document}